\documentclass{article}

\usepackage[letterpaper,nohead,top=3cm,bottom=0.5cm,left=1.5cm,right=1.5cm]{geometry}

\usepackage{amsmath,amsthm,amscd,amsfonts,amssymb,enumerate}
\usepackage{graphicx}		
\usepackage{color}
\usepackage[colorlinks]{hyperref}
\usepackage{algorithm,algorithmic}
\usepackage[titlenumbered,ruled,noend,algo2e]{algorithm2e}
\usepackage{cancel}
\usepackage{rotating}
\usepackage[graphicx]{realboxes}
\usepackage{bm}
\usepackage{tikz}
\usepackage{enumitem}
\usepackage{caption}
\usepackage{stmaryrd}
\usepackage{makecell}
\usepackage{multirow}

\usepackage[fleqn,tbtags]{mathtools}

\usepgflibrary{shapes.misc} % LATEX and plain TEX and pure pgf
\usepgflibrary[shapes.misc] % ConTEXt and pure pgf
\usetikzlibrary{shapes.misc} % LATEX and plain TEX when using Tik Z
\usetikzlibrary[shapes.misc] % ConTEXt when using Tik Z
\usetikzlibrary{shapes.geometric}

\newlength{\Oldarrayrulewidth}

\newlength{\defbaselineskip}
\newcommand{\setlinespacing}[1]%
           {\setlength{\baselineskip}{#1 \defbaselineskip}}

\renewcommand{\algorithmiccomment}[1]{\bgroup\hfill\{#1\}\egroup}

\newcommand\restr[2]{{% we make the whole thing an ordinary symbol
  \left.\kern-\nulldelimiterspace % automatically resize the bar with \right
  #1 % the function
  \littletaller % pretend it's a little taller at normal size
  \right|_{#2} % this is the delimiter
  }}
\newcommand{\littletaller}{\mathchoice{\vphantom{\big|}}{}{}{}}
\def\bkR{{\rm I\kern-.11em R}}

\theoremstyle{definition}
\newtheorem{definition}{Definition}[section]
\newtheorem{theor}{Theorem}[section]

\newcommand{\revise}[1]{#1}

\begin{document}

\newcommand{\thistitle}{Achieving $h$- and $p$-robust monolithic multigrid solvers for the Stokes equations}
\title{\thistitle}

\title{Achieving $h$- and $p$-robust monolithic multigrid solvers for the Stokes equations}

\author {Amin Rafiei\footnote{Department of Applied Mathematics, Hakim Sabzevari University, Sabzevar, Iran. E-mail address: rafiei.am@gmail.com, a.rafiei@hsu.ac.ir, arafiei@mun.ca},
Scott MacLachlan \footnote{Department of Mathematics and Statistics, Memorial University of Newfoundland, St. John's, Canada. E-mail: s.machlachlan@mun.ca},
\\
       }
%\date{}
\maketitle

\section{Introduction}
Numerical simulation of incompressible fluid dynamics has long been an important application of scientific computing. Over the years, many different discretization schemes have been considered, including the marker-and-cell finite-difference discretization \cite{MarkerAndCell}, the conforming mixed finite-element approaches of Taylor-Hood \cite{TaylorHoodINFSUP} and Scott-Vogelius \cite{ScottVogelius,ScottVogelius2} , and ${\bf H}(\text{div})$-$L^2$-conforming finite-element discretizations \cite{doi:10.1137/060649227,DiosBrezziMariniXuZikat}. Recent years has seen substantial research into these discretizations for the Stokes and Navier-Stokes equations, with a particular focus on higher-order finite-element and pressure-robust finite-element discretizations\cite{JohnMixedCGDGStokes}.  While the discretization technology has been extensively studied, there has been much less work on the development of efficient preconditioners for the linear (and linearized) systems associated with these higher-order discretizations. 

In contrast, there has been significant work on the development of block-structured preconditioners~\cite{elman, phillips2016block, wathen2020scalable} and monolithic multigrid algorithms for numerical solution of low-order discretizations of the Stokes, Navier-Stokes, and related equations.  Monolithic multigrid methods have been considered with several families of relaxation scheme, including Vanka~\cite{Vanka1,VankaType,ElementWiseExtendedVanka,VankaStokesNavierStokes}, Braess-Sarazin~\cite{BraessSarazin,ZulehnerBraessSarazin}, Uzawa~\cite{UzawaType,UzawaOlshanskii,UzawaJohnRudeWohlmuthZulehner,UzawaGmeinerHuberJohnRudeWohlmuth} and distributive~\cite{BrandtDinarDistributedRelax,DistributedRelaxBacutaVassilevskiZhang,DistributedRelaxWangChen,DistributedRelaxOosterleeGaspar} relaxation schemes. Many of these works have focused on extremely low-order discretizations, such as the marker-and-cell discretization or lowest-order Taylor-Hood discretizations. Here, we ask the question of whether such approaches can be successfully extended to higher-order finite-element discretizations, choosing to focus on monolithic multigrid with Vanka relaxation.  As is well-documented in the literature, achieving successful monolithic multigrid methods using Vanka relaxation requires careful choice of the ``Vanka patches''~\cite{JohnVankaTypeNavierStokes,ElementWiseExtendedVanka}, which have been considered for a handful of low-order discretizations, but not systematically for higher-order discretizations.

A related question is that of efficiently implementing Vanka relaxation schemes, particularly for higher-order discretizations.  While many research codes use bespoke implementations, specific to chosen problems, discretizations, and data structures, the PCPATCH implementation~\cite{PCPATCH} in PETSc was recently introduced to provide a unifying framework to implement these approaches, coupled with corresponding callbacks in Firedrake~\cite{Firedrake} to simplify patch construction for a wide variety of problems and discretizations.  Central to the PCPATCH implementation is a globally matrix-free mindset, where only the patch submatrices need to be assembled for relaxation, and this is done using callbacks to a discretization engine (Firedrake, in our case).  An alternative approach is to use a matrix-full mindset, where the global matrix is assembled once, and the patch sub-matrices are extracted from it based on lists of degree of freedom (DoF) indices.  This is implemented in Firedrake's ASMPatchPC class, which we will compare in performance to the PCPATCH implementation.

Another key motivation for considering a matrix-full implementation is that it is more easily extensible than the callback-driven implementation of PCPATCH, which requires information about the mesh elements associated with and surrounding the patch in order to correctly assemble the local system.  Since this work has already been done in assembling the global matrix, no additional work (other than identifying the DoFs to include in each patch) is needed for matrix-full implementations of patch-based relaxation schemes.  In particular, this lets us implement and generalize ``extended'' Vanka patches~\cite{ElementWiseExtendedVanka} for higher-order ${\bf H}(\text{div})$-$L^2$-conforming discretizations that are not currently possible in PCPATCH.  At the same time, we show that both ASMPatchPC and PCPATCH can be easily leveraged to implement ``composite'' patches that prove effective for higher-order Taylor-Hood discretizations on both triangular and quadrilateral meshes.

While the numerical results that follow show robust iteration counts and times-to-solution for many of the discretizations when using monolithic multigrid as a preconditioner for FGMRES, a key open question is that of reliable stopping criteria for iterative methods such as these.  This question is not unique to the Stokes equations, and has recently been considered also for higher-order discretizations of the Poisson problem~\cite{guo2023stopping}, where several possible solutions are considered.  Here, we note only that standard residual-based stopping tolerances prove ineffective for higher-order discretizations of the Stokes equations on fine-enough meshes, leading to the need for future work to extend the results of Guo et al.~\cite{guo2023stopping} and similar research to other settings.

The remainder of this paper is organized as follows. In Section~\ref{sec:mixedFEM}, we review the Stokes equations and their discretization using mixed finite-element methods, paying particular attention to the distinction between pressure-robust and other discretizations.  In Section~\ref{sec:monolithicMG}, we propose monolithic geometric multigrid solvers for higher-order stable conforming and nonconforming discretizations of the Stokes equations. The implementation of ASMPatchPC-based Vanka relaxation is a key ingredient for constructing these algorithms, which is reviewed in this section. In Section~\ref{sec:numerics}, we show performance of this monolithic multigrid framework as preconditioners for FGMRES, verifying the robustness of these solvers with respect to iteration counts and time-to-solution, and giving insight into the choice of stopping criteria for these solvers. Conclusions and future work are discussed in Section~\ref{sec:conclusions}.

%%%%%%%%%%%%%%%%%%%%%%%
\section{Mixed finite-element discretization of Stokes equations}\label{sec:mixedFEM}
The Stokes equations describe steady, incompressible viscous flows. For a simply connected polygonal domain $\Omega \subseteq \mathbb{R}^2$, the (time-steady) Stokes equations are given by
\begin{eqnarray}\label{eq:Stokes}
-\nabla \cdot (2\nu \bm{\varepsilon}({\bf u})) + \nabla p & = & {\bf f}\text{ in }\Omega, \label{eq:Stokes1}\\
\nabla\cdot{\bf u} & = & 0\text{ in }\Omega, \label{eq:StokesDivFree}
\end{eqnarray}
\noindent where $\nu$ is the kinematic fluid viscosity, ${\bf u}:\Omega\rightarrow \mathbb{R}^2$ is the velocity field, $p: \Omega \rightarrow \mathbb{R}$ is the pressure field, and  ${\bf f}:\Omega\rightarrow \mathbb{R}^2$ is the external body force acting on the fluid. Here, $\bm{\varepsilon}({\bf u}) = \frac{1}{2}(\nabla {\bf u} + \nabla {\bf u}^T)$ is the symmetric strain-rate tensor, where $\nabla {\bf u}$ is the tensor defined by $(\nabla {\bf u})_{ij} = \frac{\partial u_i}{\partial x_j}$. Continuum equations~\eqref{eq:Stokes1} and~\eqref{eq:StokesDivFree} are referred to as the momentum and continuity balance equations, respectively, with~\eqref{eq:StokesDivFree} representing the incompressibility of the fluid.

No-slip (essential) Dirichlet boundary conditions are formulated as
\begin{eqnarray}\label{eq:StokesBCforCGFEM}
{\bf u} & = & {\bf 0}\text{ on }\partial \Omega.  \label{eq:StokesDirichletNoSlip}
\end{eqnarray}
In the case of enclosed flow, an alternative to~\eqref{eq:StokesBCforCGFEM} is the no-flux condition,  
\begin{eqnarray}\label{eq:StokesBCforDGFEM1}
{\bf u\cdot n} = {\bf 0}\text{ on }\partial \Omega, \label{eq:StokesDirichletEnclosedNoFlux}
\end{eqnarray}
where ${\bf n}$ is the outward unit normal vector to the boundary. This essential boundary condition can be complemented by imposing the natural boundary condition  
\begin{eqnarray}\label{eq:StokesBCforDGFEM2}
( ( 2\nu {\bm \varepsilon}({\bf u}) -p {\bf I} ) {\bf n} )\cdot {\bf t} = {\bf 0},\text{ on }\partial \Omega, \label{eq:StokesDirichletTangentNormalStress}
\end{eqnarray}
on the tangential component of the normal stresses, where ${\bf t}$ is the tangential unit vector on the boundary and ${\bf I}$ is the identity tensor. Since ${\bf n}\cdot{\bf t} = 0$, \eqref{eq:StokesBCforDGFEM2} reduces to $({\bm \varepsilon}({\bf u}){\bf n})\cdot {\bf t} = 0$.

The standard weak formulation of~\eqref{eq:Stokes1} and~\eqref{eq:StokesDivFree} relies on the spaces
\begin{eqnarray*}
  {\bf H_0^1}(\Omega) & = & \left\{ {\bf v}\in {\bf H^1}(\Omega)~\middle|~{\bf v} = {\bf 0}\text{ on }\partial \Omega\right\}, \\
  L_0^2(\Omega) & = & L^2(\Omega) / \mathbb{R} := \left\{q\in L^2(\Omega)~\middle|~ \int_{\Omega}qdx = 0\right\},
\end{eqnarray*}
with their associated norms $\|\cdot\|_{{\bf H^1}(\Omega)}$ and $\|\cdot\|_{L^2(\Omega)}$, where we use ${\bf H^1}(\Omega) = [H^1(\Omega)]^2$ to denote vector-valued functions with components in $H^1(\Omega)$.  For a given ${\bf f}\in {\bf H^{-1}}(\Omega)$, the variational formulation is to find $({\bf u},p)  \in {\bf H_0^1}(\Omega)\times L_0^2(\Omega)$ such that
\begin{eqnarray}\label{eq:VariationalStokes1}
a({\bf u},{\bf v}) + b({\bf v},p) & = & \langle{\bf f},{\bf v}\rangle, \nonumber\\
           b({\bf u},q) & = & 0, 
\end{eqnarray}
for all ${\bf v}\in {\bf H_0^1}(\Omega)$ and $q\in L_0^2(\Omega)$, with bilinear forms
\begin{eqnarray}\label{eq:VariationalStokesA}
a({\bf u},{\bf v}) & := & 2\nu \int_{\Omega}\bm{\varepsilon}({\bf u}): \bm{\varepsilon}({\bf v})~dx, \label{eq:VariationalStokesa}\\
b({\bf v},p) & := & - \int_{\Omega} (\nabla\cdot {\bf v}) p~dx, \label{eq:VariationalStokesb}
\end{eqnarray}
and linear form
\[
\langle{\bf f},{\bf v}\rangle := \int_{\Omega} {\bf f}\cdot {\bf v}~dx.
\]
These forms are well-known to be continuous, and $a({\bf u},{\bf v})$ is coercive on ${\bf H_0^1}(\Omega)$.  To prove uniqueness of solutions, we also need that $b({\bf v},p)$ satisfies an inf-sup condition, which is again well-known~\cite[Example 4.2.2]{MixedFEMAppl}.

At the discrete level, we let $\mathcal{T}_h$ be a triangulation of $\Omega$, and first consider
${\bf V}_h = {\bf V}_h(\mathcal{T}_h)$ and $Q_h = Q_h(\mathcal{T}_h)$ to be ${\bf H_0^1}$- and $L_0^2$-conforming finite-element spaces over $\mathcal{T}_h$, emphasizing that all components of vector-valued functions ${\bf v}_h\in {\bf V}_h$ are continuous across all inner-element boundaries of $\mathcal{T}_h$.  Using conforming discretizations has the advantage that the discrete weak forms can be kept the same as continuous ones, leading to the discrete form of~\eqref{eq:VariationalStokes1} as finding approximate solution $({\bf u}_h,p_h)\in ({\bf V}_h,Q_h)$ such that
\begin{eqnarray}\label{eq:VariationalDiscStokes}
a({\bf u}_h,{\bf v}_h) + b({\bf v}_h,p_h) & = & \langle{\bf f},{\bf v}_h\rangle, \label{eq:DiscMomentum} \nonumber\\
           b({\bf u}_h,q_h) & = & 0, \label{eq:DiscCompressibility}
\end{eqnarray}
for all ${\bf v}_h\in {\bf V}_h$ and $q_h\in Q_h$.  Since the continuity of the bilinear forms and coercivity of $a({\bf u}_h,{\bf v}_h)$ are inherited from the continuum, we only need to check the inf-sup condition at the discrete level, leading to the following existence and uniqueness result.
%%%%%%%%%%%%%%%%%%%%%%
\begin{theor} \label{Th:PressRobust}
Let $({\bf u},p)\in {\bf H_0^1}(\Omega)\times L_0^2(\Omega)$ be the solution of \eqref{eq:VariationalStokes1} and suppose that
\begin{eqnarray}\label{eq:DiscINFSUP}
\adjustlimits\inf_{q_h\in Q_h\backslash\{0\}} \sup_{{\bf v}_h\in {\bf V}_h \backslash\{\bf 0\}} \frac{b_h({\bf v}_h,q_h)}{\|{\bf v}_h\|_{{\bf H^1}(\Omega)}\|q_h\|_{L^2(\Omega)}}\ge \beta_h>0.
\end{eqnarray}
Then, there exists a unique solution $({\bf u}_h,p_h)\in ({\bf V}_h,Q_h)$ of~\eqref{eq:VariationalDiscStokes}, and the following error bounds are satisfied 
\begin{eqnarray}
\|{\bf u}_h - {\bf u}\|_{{\bf H^1}(\Omega)} & \le &  \left(\frac{2\|a\|}{\gamma} + \frac{2\|a\|^{1/2}\|b\|}{(\gamma)^{1/2}\beta_h}  \right)E_{\bf u} + \left(\frac{\|b\|}{\gamma}\right)E_p, \label{eq:BoundExErrVelocit}\\
\|p_h - p\|_{L^2(\Omega)} & \le &  \left(\frac{2\|a\|^{3/2}}{(\gamma)^{1/2}\beta_h} + \frac{\|a\|\|b\|}{\beta_h^2} \right)E_{\bf u} + \left(\frac{3\|a\|^{1/2}\|b\|}{(\gamma)^{1/2}\beta_h}\right)E_p, \label{eq:BoundExErrPressure}
\end{eqnarray}
where $\gamma$ is the coercivity constant of $a({\bf u},{\bf v})$ and the approximation errors, $E_{\bf u}$ and $E_p$, are defined as
\begin{eqnarray*}
E_{\bf u} & := & \inf\limits_{{\bf v}_h\in {\bf V}_h} \|{\bf u} - {\bf v}_h\|_{{\bf H^1}(\Omega)}, \nonumber \\
E_{p} & := & \inf\limits_{q_h\in Q_h} \|p - q_h\|_{L^2(\Omega)}. 
\end{eqnarray*}
{\bf Proof:} See \cite[Chapter 8]{MixedFEMAppl}. 
\end{theor}
%%%%%%%%%%%%%%%%%%%%%%

Let $P_k = P_k(\mathcal{T}_h)$ be the continuous finite-element space on triangles with piecewise polynomial approximation of degree $k$ and ${\bf P_k}$ be its vector-valued counterpart.  Similarly, let $Q_k = Q_k(\mathcal{T}_h)$ be the continuous finite-element space on quadrilaterals with piecewise bi-polynomial approximation of degree $k$ and ${\bf Q_k}$ be its vector-valued counterpart.  Among the most common conforming finite-element schemes for~\eqref{eq:DiscMomentum} are the Taylor-Hood finite-element pairs, $({\bf V}_h,Q_h) = ({\bf Q_k},Q_{k-1})$ (on quadrilaterals) and $({\bf V}_h,Q_h) = ({\bf P_k},P_{k-1})$ (on triangles), for $k\ge 2$.  For any fixed $k\geq 2$, these elements are known to satisfy an $h$-uniform version of the discrete inf-sup condition~\eqref{eq:DiscINFSUP} with $\beta_h\ge \beta_0 > 0$ for all $h>0$, on any mesh~\cite[Section 8.8.2]{MixedFEMAppl}. Since interpolation error in the ${\bf H^1}$ norm for ${\bf P_k}$ or ${\bf Q_k}$ is of order $h^k$ (where $h$ is, for example, the maximum diameter of any element in $\mathcal{T}_h$) and interpolation error in the $L^2$ norm for $P_{k-1}$ or $Q_{k-1}$ is also of order $h^k$, \eqref{eq:BoundExErrVelocit} and~\eqref{eq:BoundExErrPressure} lead to order $h^k$ convergence for both velocity and pressure (in their respective natural norms), assuming smooth-enough solutions to the Stokes equations.

In~\eqref{eq:BoundExErrVelocit}, the velocity approximation error is bounded by interpolation errors for both the velocity and pressure.  While this seems natural from a finite-element perspective, it sits in contrast to a fundamental invariance property of the continuous equations~\eqref{eq:Stokes1} and \eqref{eq:StokesDivFree}, in which incrementing the external force by a gradient field only changes the pressure solution and not the velocity field, with
\begin{eqnarray*}
{\bf f} \rightarrow {\bf f} + \nabla \psi \implies ({\bf u},p) \rightarrow ({\bf u},p+\psi).
\end{eqnarray*}
This invariance property is also inherited by the discretized equation form in~\eqref{eq:DiscMomentum}, so long as $\psi\in P_{k-1}$ or $Q_{k-1}$, but the upper bound in inequality~\eqref{eq:BoundExErrVelocit} does not reflect this invariance property. Thus, there has been significant recent research to improve this discretization framework, through the introduction of ``pressure-robust'' schemes~\cite{JohnMixedCGDGStokes}.

\begin{definition}\label{def:PressRobust}
Consider the mixed finite-element dicretization of the Stokes equations in~\eqref{eq:DiscCompressibility}.
If the divergence constraint in~\eqref{eq:StokesDivFree} is enforced exactly and the velocity error does not depend on the pressure, then such method is said to be a {\bf pressure-robust} method. 
\end{definition}

There are two important consequences of using a pressure-robust method, both to remove the dependence on the pressure interpolation error in the bound on the velocity approximation error, and to result in a discrete velocity field, ${\bf u}_h$, that satisfies the continuous incompressibility constraint~\eqref{eq:StokesDivFree} (meaning the velocity is strongly and pointwise divergence-free).  The natural way to achieve a pressure-robust method is to ensure that $\nabla\cdot{\bf u}_h \in Q_h$ (the pressure space) for any ${\bf u}_h \in {\bf V}_h$.  It is easy to see this is not the case when we consider Taylor-Hood elements, since the divergence of a continuous piecewise-polynomial vector field, ${\bf u}_h$, need not be continuous.  Thus, we conclude that the Taylor-Hood elements are not pressure-robust, and note that the quantity $\|\nabla\cdot{\bf u}_h\|_{L^2(\Omega)}$ can, indeed, become quite large~\cite{CaseErvinLinkeRebholtz}. 

One option to achieve pressure-robust discretization while keeping the velocity in a standard space is to change the pressure space to discontinuous piecewise polynomial functions, defining $dP_k = dP_k(\mathcal{T}_h)$ for triangular meshes.  This leads to the Scott-Vogelius spaces on simplices, $({\bf V}_h,Q_h) = ({\bf P_k},dP_{k-1})$~\cite{ScottVogelius, ScottVogelius2}.  These spaces, however, are only known to be stable under assumptions on the mesh; for example, in two dimensions, the triangular elements are stable for $k\geq 2$ on triangular meshes that result from a single step of barycentric refinement of any triangular mesh~\cite{qin1994convergence}.  Higher-order pairs are known to be stable under weaker assumptions.  In two dimensions, these elements are stable for $k\geq 4$ if the mesh does not have any nearly singular vertices~\cite{ScottVogelius}.  While these elements are attractive, since they offer the same convergence orders as Taylor-Hood, the restrictions on meshes for which they are inf-sup stable and difficulties in defining effective preconditioners on those meshes~\cite{farrell2021reynolds} makes direct comparison between solvers for Taylor-Hood and Scott-Vogelius somewhat complicated.  For this reason, we do not consider Scott-Vogelius elements further in this work.

%%%%%%%%%%%%%%%%%%%%%%%%%%%%%%%%%%%%
\subsection{Stable nonconforming discretizations}
%%%%%%%%%%%%%%%%%%%%%%%%%%%%%%%%%%%%

To achieve inf-sup stable discretizations on general meshes, we move to the nonconforming case, where ${\bf V}_h \not\subset {\bf H_0^1}$.  In particular, we will focus on the class of ${\bf H}(\text{div})$-$L^2$-conforming discretizations~\cite{doi:10.1137/060649227, DiosBrezziMariniXuZikat}. The Sobolev space ${\bf H}(\text{div},\Omega)$ consists of vector fields for which the components of the field and its weak divergence
are square-integrable over $\Omega$. This space and its corresponding norm are defined as
\begin{eqnarray*}
{\bf H}(\text{div},\Omega) := \left\{ {\bf v}\in {\bf L^2}(\Omega)~\middle|~\nabla \cdot {\bf v}\in L^2(\Omega) \right\}, & & \|{\bf v}\|_{{\bf H}(\text{div},\Omega)}^2 := \|{
\bf v}\|^2_{{\bf L^2}(\Omega)} + \|\nabla \cdot {\bf v}\|^2_{L^2(\Omega)}.
\end{eqnarray*}
As above, we will consider the subspace of this space that applies relevant boundary conditions, with
\begin{eqnarray*}
{\bf H_0}(\text{div},\Omega) := \left\{ {\bf v}\in {\bf H}(\text{div},\Omega)~\middle|~{\bf v} \cdot {\bf n} = 0\text{ on }\partial \Omega \right\}.
\end{eqnarray*}

Given a triangulation, $\mathcal{T}_h$, of the domain $\Omega$, two standard ${\bf H}(\text{div})$-conforming finite-element spaces on $\mathcal{T}_h$ are the Raviart-Thomas~\cite{RTelement} and Brezzi-Douglas-Marini~\cite{BDMelements} spaces.  Both spaces consist of piecewise polynomial basis functions with continuous normal components across the edges of $\mathcal{T}_h$, but differ in how the polynomial spaces on each triangle are defined.  For a given triangle, $T\in\mathcal{T}_h$, we define
\[
\mathcal{RT}_k(T) = \left\{{\bf v}\in L^2(T)~\middle|~{\bf v} \in {\bf P_{k-1}}(T) + {\bf x}P_{k-1}(T)\right\} \text{ and } \mathcal{BDM}_k(T) = \left\{{\bf v}\in L^2(T)~\middle|~{\bf v} \in {\bf P_{k}}(T)\right\},
\]
where ${\bf x}$ is the standard coordinate vector on $T$.
We define $\mathcal{RT}_k = \mathcal{RT}_k(\mathcal{T}_h)$ and $\mathcal{BDM}_k = \mathcal{BDM}_k(\mathcal{T}_h)$ to be the global Raviart-Thomas and Brezzi-Douglas-Maraini finite-element spaces over $\mathcal{T}_h$, respectively. 
These are constructed so that $\mathcal{RT}_k$ and $\mathcal{BDM}_k$ are ${\bf H}(\text{div})$-conforming spaces, in the sense that a vector-valued function ${\bf v}\in \mathcal{RT}_k$ or ${\bf v}\in \mathcal{BDM}_k$ need not be continuous in all components over $\mathcal{T}_h$, but its normal component must be over each edge between two triangles in $\mathcal{T}_h$.  Strongly imposing Dirichlet boundary conditions on ${\bf v}\cdot{\bf n}$ for functions in these spaces is made possible by suitable selection of the degrees of freedom, including dual basis functions corresponding to point or integral moment evaluation of the normal component along edges of the mesh.  When we strongly impose ${\bf v}\cdot{\bf n} = 0$ on $\partial\Omega$, we will refer to the space as an ${\bf H_0}(\text{div})$-conforming space.

In developing a variational formulation for solving the Stokes equations \eqref{eq:Stokes1} and \eqref{eq:StokesDivFree} with boundary conditions \eqref{eq:StokesBCforDGFEM1} and \eqref{eq:StokesBCforDGFEM2}, one possibility is to search for solutions to a variational formulation in a space of divergence-free functions, thereby satisfying \eqref{eq:StokesDivFree} by construction. However, this does not translate well to finite-element formulations.  Instead, we make use of standard interior penalty techniques to define a suitable weak formulation.

For $s\ge 1$, we can define the ``broken'' Sobolev spaces on triangulation $\mathcal{T}_h$ as
\begin{eqnarray*}
H^s_b(\mathcal{T}_h) = \left\{\phi\in L^2(\Omega)~\middle|~\restr{\phi}{T}\in H^s(T),\quad \forall T\in \mathcal{T}_h\right\},
\end{eqnarray*}
denoting the vector and tensor analogues of $H^s_b(\mathcal{T}_h)$ by ${\bf H^s_b}(\mathcal{T}_h)$ and $\bm{\mathcal{H}}^{\bf s}_{\bf b}(\mathcal{T}_h)$, respectively. 
Let $\xi_h^{I}$ be the set of internal edges of $\mathcal{T}_h$, and let $e\in \xi_h^{I}$ be adjacent to two elements, $T_1$ and $T_2$. Let ${\bf n_1}$ and ${\bf n_2}$ denote the outward unit normal vectors on $e$ associated with $T_1$ and $T_2$. For a tensor field $\bm{\tau}\in \bm{\mathcal{H}^{\bf 1}_{\bf b}}({\mathcal{T}_h})$, the average operator on $e$ is defined as 
\begin{eqnarray*}
\{\bm{\mathcal{\tau}}\}_e = \frac{\bm{\tau}_1 + \bm{\tau}_2}{2},
\end{eqnarray*}
where $\bm{\tau}_i = \restr{\bm{\tau}}{\partial T_i}$ is the trace of $\bm{\tau}$ on $e$ from triangle $T_i$. For a vector field ${\bf v}\in {\bf H^1_b}(\mathcal{T}_h)$, the symmetric (matrix-valued) jump operator on $e\in \xi_h^{I}$ is given by
\begin{eqnarray*}
\llbracket{\bf v}\rrbracket_e = {\bf v_1}\odot {\bf n_1} + {\bf v_2}\odot {\bf n_2},
\end{eqnarray*}
where ${\bf v}\odot {\bf n} = \frac{{\bf v}{\bf n^T} + {\bf n}{\bf v^T}}{2}$ and, as above, ${\bf v}_i$ is the trace of ${\bf v}$ on $e$ from triangle $T_i$. 

Let ${\bf u}_h$, ${\bf v}_h \in {\bf V}_h$ for some ${\bf H_0}(\text{div})$-conforming finite-element space, ${\bf V}_h$, on some mesh $\mathcal{T}_h$, and note that this implies that both $\bm{\varepsilon}({\bf u}_h) \in \bm{\mathcal{H}}^{\bf 1}_{\bf b}(\mathcal{T}_h)$ and ${\bf v}_h \in {\bf H^1_b}(\mathcal{T}_h)$.  Working elementwise, we can integrate by parts to get that
\begin{equation*}
\sum_{T\in \mathcal{T}_h} -\int_{T} \left[\nabla \cdot (2\nu \bm{\varepsilon}({\bf u}_h))\right]\cdot {\bf v}_h~dx = \sum_{T\in \mathcal{T}_h} 2\nu \int_{T}\bm{\varepsilon}({\bf u}_h): \bm{\varepsilon}({\bf v}_h)~dx - \sum_{T\in \mathcal{T}_h } 2\nu \int_{\partial T} {\bf v}_h\cdot\left[\bm{\varepsilon}({\bf u}_h) {\bf n}_T\right]~ds,
\end{equation*}
where ${\bf n}_T$ is the outward unit normal vector on $\partial T$.  A standard identity is, then, that
\begin{equation*}
  \sum_{T\in \mathcal{T}_h } 2\nu \int_{\partial T} {\bf v}_h\cdot\left[\bm{\varepsilon}({\bf u}_h) {\bf n}_T\right]~ds = \sum_{e\in \xi_h^{I}} 2\nu \int_{e} \llbracket{\bf v}_h\rrbracket_e:\{\bm{\varepsilon}({\bf u}_h)\}_e~ds,
\end{equation*}
where we use the boundary conditions in~\eqref{eq:StokesBCforDGFEM1} and~\eqref{eq:StokesBCforDGFEM2} to eliminate contributions from edges on $\partial\Omega$ in the summation on the right.
Noting, further, that ${\bf v}_h \in {\bf V}_h$ implies that the normal component of ${\bf v}_h$ must be continuous across all $e\in \xi_h^{I}$, allows us to simplify the term on the right-hand side to just include the tangential component of ${\bf v}_h$, which we denote as $\llbracket{\bf v}_h\rrbracket_e^{\bf t} = \llbracket\left({\bf v}_h\cdot{\bf t}\right){\bf t}\rrbracket_e$, where ${\bf t}$ is the unit tangent vector to $e$.  This gives us
\begin{equation*}
\sum_{T\in \mathcal{T}_h} -\int_{T} \left[\nabla \cdot (2\nu \bm{\varepsilon}({\bf u}_h))\right]\cdot {\bf v}_h~dx = \sum_{T\in \mathcal{T}_h} 2\nu \int_{T}\bm{\varepsilon}({\bf u}_h): \bm{\varepsilon}({\bf v}_h)~dx -\sum_{e\in \xi_h^{I}} 2\nu \int_{e} \llbracket{\bf v}_h\rrbracket_e^{\bf t}:\{\bm{\varepsilon}({\bf u}_h)\}_e~ds.
\end{equation*}
To complete the interior penalty discretization, we add a consistent symmetrization term,
\begin{equation*}
-\sum_{e\in \xi_h^{I}} 2\nu \int_{e} \llbracket{\bf u}_h\rrbracket_e^{\bf t}:\{\bm{\varepsilon}({\bf v}_h)\}_e~ds,
\end{equation*}
and a suitable penalty term,
\begin{equation*}
2\nu \alpha  \sum_{e\in \xi_h^{I}} h_e^{-1}\int_{e} \llbracket {\bf u}_h \rrbracket_e^{\bf t}:\llbracket {\bf v}_h \rrbracket_e^{\bf t}~ds,
\end{equation*}
where $\alpha$ is a penalty parameter and $h_e$ is the length of edge $e$, both justified by noting that $\llbracket {\bf u}_h \rrbracket_e^{\bf t}$ is zero for a solution ${\bf u}_h \in {\bf H^1}(\Omega)$, but not necessarily for a nonconforming solution ${\bf u}_h \in {\bf V}_h$.  This leads to the symmetric bilinear form
\begin{eqnarray}\label{eq:DiscontBilinearForm}
{a}_h({\bf u}_h,{\bf v}_h) & := & 2\nu \left[  \sum_{T\in \mathcal{T}_h} \int_{T} \bm{\varepsilon}({\bf u}_h):\bm{\varepsilon}({\bf v}_h)~dx - \sum_{e\in \xi_h^{I}}\int_{e} \{ \bm{\varepsilon}({\bf u}_h) \}_e : 
\llbracket{\bf v}_h\rrbracket_e^{\bf t}~ds - \sum_{e\in \xi_h^{I}}\int_{e} \llbracket{\bf u}_h\rrbracket_e^{\bf t}:\{\bm{\varepsilon}({\bf v}_h)\}_e~ds \right] \nonumber\\
                                   & + & 2\nu \alpha  \sum_{e\in \xi_h^{I}} h_e^{-1}\int_{e} \llbracket {\bf u}_h\rrbracket_e^{\bf t}:
\llbracket{\bf v}_h\rrbracket_e^{\bf t}~ds.
\end{eqnarray}

With~\eqref{eq:DiscontBilinearForm}, we can write the weak form for~\eqref{eq:Stokes1} and~\eqref{eq:StokesDivFree} for ${\bf H_0}(\text{div})$-conforming space ${\bf V}_h$ and $L^2_0$-conforming space $Q_h$ as finding $({\bf u}_h,p_h)\in ({\bf V}_h,Q_h)$ such that
\begin{eqnarray}\label{eq:VariationalDiscStokesHdivConform}
{a}_h({\bf u}_h,{\bf v}_h) + b({\bf v}_h,p_h) & = & \langle{\bf f},{\bf v}_h\rangle,\quad \forall {\bf v}_h\in {\bf V}_h, \label{eq:DiscMomentumHdivConform} \nonumber\\
           b({\bf u}_h,q_h) & = & 0,\quad \forall q_h\in Q_h, \label{eq:DiscCompressibilityHdivConform}
\end{eqnarray}
where $b(\cdot,\cdot)$ is defined as in \eqref{eq:VariationalStokesb}.  
Defining
\begin{equation*}
\|{\bf v}\|_{DG}^2 = 2\nu \left(\sum_{T\in \mathcal{T}_h} \|\nabla {\bf v}\|_{{\bf L^2}(T)}^2 + \sum_{e\in \xi_h^I} h_e^{-1}\int_e \llbracket {\bf v}\rrbracket_e^{\bf t}:\llbracket{\bf v}\rrbracket_e^{\bf t}~ds\right),
\end{equation*}
we can show that this is a seminorm on ${\bf H^1_b}(\mathcal{T}_h)$ and a norm on ${\bf H_0}(\text{div},\Omega)$.
To address the difference in  boundary conditions on ${\bf H_0}(\text{div})$-conforming spaces versus those on ${\bf H_0^1}$-conforming spaces, we define
\begin{equation*}
{\bf H^1_{0,n}}(\Omega) = \left\{ {\bf v}\in {\bf H^1}(\Omega)~\middle|~{\bf v}\cdot{\bf n} = 0 \text{ on }\partial\Omega\right\},
\end{equation*}
and consider the continuous saddle-point problem in~\eqref{eq:VariationalStokes1} with ${\bf u},{\bf v}$ restricted to ${\bf H_{0,n}^1}(\Omega)$, with the natural boundary condition in~\eqref{eq:StokesBCforDGFEM2} also imposed.  The saddle-point problem restricted to ${\bf H_{0,n}^1}(\Omega)$ has a unique solution by standard analytical techniques~\cite{DiosBrezziMariniXuZikat}.
The following results are known in the literature~\cite{DiosBrezziMariniXuZikat}.
\begin{theor} \label{Th:CoercContinDiscStokesHdivConform}
Let $({{\bf u}},{p})  \in {\bf H_{0,n}^1}(\Omega)\times L_0^2(\Omega)$ be the solution of the continuous saddle-point problem \eqref{eq:VariationalStokes1}, with ${\bf u},{\bf v}$ restricted to ${\bf H_{0,n}^1}(\Omega)$. 
Also let $({\bf V}_h,Q_h)$ be an ${\bf H_0}(\text{div})$-$L^2_0$-conforming finite-element pair equipped with the norms $\|\cdot\|_{DG}$ and $\|\cdot\|_{L^2(\Omega)}$, respectively. Then, for sufficiently large penalty parameter, $\alpha$,
\begin{itemize}
\item the bilinear form ${a}_h(\cdot,\cdot)$ is continuous with respect to the $\|\cdot\|_{DG}$ norm and $b(\cdot,\cdot)$ is continuous with respect to the $\|\cdot\|_{DG}$ and $\|\cdot\|_{L^2(\Omega)}$ norms,
\item there is a constant $\gamma>0$ such that the coercivity condition 
\begin{eqnarray*}
{a}_h({\bf v}_h,{\bf v}_h) \ge \gamma \|{\bf v}_h\|_{DG}^2, & \forall {\bf v}_h\in {\bf V}_h,
\end{eqnarray*}
is satisfied, and 
\item there exits a constant $\beta>0$, independent of the mesh size $h$, such that the inf-sup condition 
\begin{eqnarray*}
\adjustlimits\inf_{q_h\in Q_h\backslash\{0\}} \sup_{{\bf v}_h\in {\bf V}_h \backslash\{\bf 0\}} \frac{b({\bf v}_h,q_h)}{\|{\bf v}_h\|_{DG}~~\|q_h\|_{L^2(\Omega)}}\ge \beta>0, 
\end{eqnarray*}
is satisfied.
\end{itemize}
Thus, the solution $({\bf u}_h,p_h)$ to~\eqref{eq:VariationalDiscStokesHdivConform} exists and is unique.  Furthermore,
\begin{eqnarray*}
\nabla\cdot {\bf u}_h = 0,\text{ in }\Omega,
\end{eqnarray*}
and there exists a constant $C$, independent of $h$, such that for every ${\bf v}_h\in {\bf V}_h$ with $\nabla\cdot {\bf v}_h = 0$ and for every $q_h\in Q_h$, the following error estimates hold,
\begin{eqnarray*}
\|{{\bf u}} - {\bf u}_h\|_{DG} & \le & C\|{{\bf u}} - {\bf v}_h \|_{DG}, \nonumber\\
\|{p} - p_h\|_{L^2(\Omega)} & \le & C\left( \|{p} - q_h\|_{L^2(\Omega)} + \| {{\bf u}} - {\bf v}_h\|_{DG}\right).
\end{eqnarray*}
\end{theor}

%%%%%%%%%%%%%%%%%%%%%%%
\section{Monolithic multigrid solvers}\label{sec:monolithicMG}
In the numerical results that follow, we consider four discretizations of the Stokes equations \eqref{eq:Stokes1} and \eqref{eq:StokesDivFree}.  Two of these are conforming Taylor-Hood discretizations,
\begin{itemize}
\item {\bf Case 1}: solving~\eqref{eq:VariationalDiscStokes} on a triangular mesh with $({\bf V}_h,Q_h) = ({\bf P_k},P_{k-1})$, for $k\ge 2$,
\item {\bf Case 2}:  solving~\eqref{eq:VariationalDiscStokes} on a quadrilateral mesh with $({\bf V}_h,Q_h) = ({\bf Q_k},Q_{k-1})$, for $k\ge 2$.
\end{itemize}
Two are non-conforming ${\bf H}(\text{div})$-$L^2$ discretizations,
\begin{itemize}
\item {\bf Case 3}: solving~\eqref{eq:VariationalDiscStokesHdivConform} on a triangular mesh with $({\bf V}_h,Q_h) = (\mathcal{BDM}_{k},{dP}_{k-1})$, for $k\ge 1$,
\item {\bf Case 4}:  solving~\eqref{eq:VariationalDiscStokesHdivConform} on a triangular mesh with $({\bf V}_h,Q_h) = (\mathcal{RT}_{k},{dP}_{k-1})$, for $k\ge 1$.
\end{itemize}
For each case, we consider solution of the global linear system $\mathcal{A}x = b$ using a monolithic $h \mkern 2mu$-multigrid $V$-cycle as a preconditioner for FGMRES.  We note that we use FGMRES not because of its flexibility (as we use stationary preconditioners), but because of its explicit storage of the preconditioned Arnoldi vectors.  While this requires extra vector storage, it lowers the cost of the solution process in comparison to classical right-preconditioned GMRES, which requires an extra application of the preconditioner once the solution over the Krylov space is found~\cite{YousefBook}.

\begin{figure}[t]
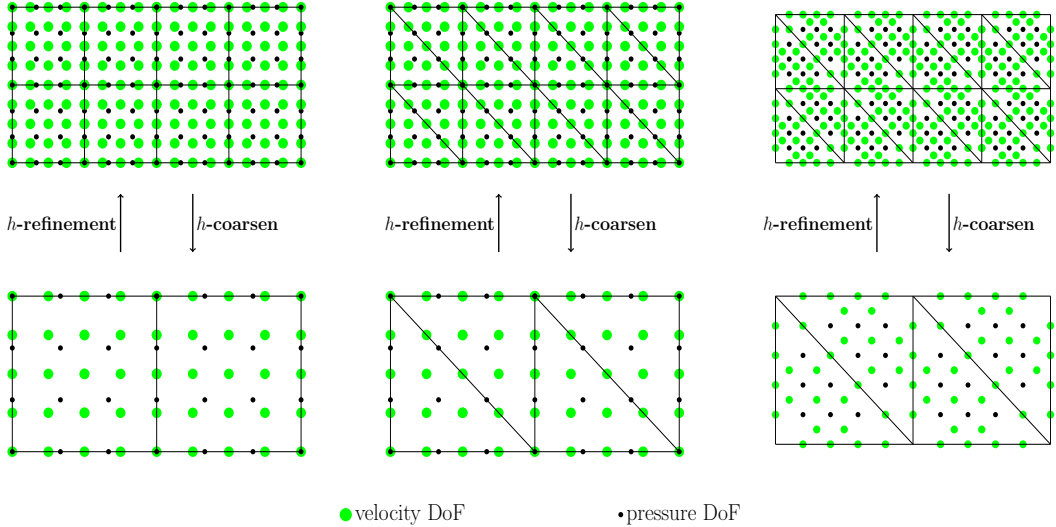
 
\begin{center}

\end{center}
\caption{$h \mkern 2mu$-multigrid coarsening and refinement for $({\bf Q_4},Q_3)$ ({\bf left}), $({\bf P_4},P_3)$ ({\bf center}) and $(\mathcal{BDM}_3,{dP_2})$ ({\bf right}) discretizations of the Stokes equations.
Big (green) and small (black) circles are associated with velocity and pressure DoFs, respectively.}
\label{fig:GeomHRefineStokes}
\end{figure}

Here, we review multigrid as a solver for the global linear system.
For a fixed value of polynomial order, $k$, let $({\bf V}_{h,1},Q_{h,1})$ be finite-element spaces of that order defined on mesh $\mathcal{T}_1$. Based on standard geometric $h$-refinement, we can build a quasi-uniform family of nested grids, $\mathcal{T}_1\subset \mathcal{T}_2\subset \cdots \subset \mathcal{T}_L$, for domain $\Omega$, where $\mathcal{T}_L$ is the finest mesh.
The nested grids $\mathcal{T}_1$, $\mathcal{T}_2$, $\cdots$, $\mathcal{T}_L$ induce the nested pairs of finite-element spaces 
\begin{eqnarray*}
({\bf V}_{h,\ell -1},Q_{h,\ell -1}) \subset ({\bf V}_{h,\ell },Q_{h,\ell}) & & \ell = 2, \cdots, L,
\end{eqnarray*}
where $({\bf V}_{h,\ell },Q_{h,\ell})$ are the discretization spaces for velocity and pressure on  mesh $\mathcal{T}_{\ell}$.
 This technique is depicted in Figure \ref{fig:GeomHRefineStokes} for meshes associated with $({\bf Q_4},Q_3)$, $({\bf P}_4,P_3)$, and $(\mathcal{BDM}_3,{dP}_2)$.

%%%%%%%%%%%%%%%%%%%%%%%%%%%%%%%
%%%%%%%%%%%%%%%%%%%%%%%%%%%%%%%%%%%%%

Discretization of the weak variational forms in \eqref{eq:VariationalDiscStokes} and \eqref{eq:VariationalDiscStokesHdivConform} at level $\ell$ of the multigrid  hierarchy, using finite-element spaces ${\bf V}_{h,\ell}$ and $Q_{h,\ell}$, produces a linear system of saddle-point equations that we write as
\begin{eqnarray*}
\underbrace{\left( \begin{array}{cc}
A_{\ell} & B_{\ell}^T \\
B_{\ell} & 0
\end{array}\right)}_{\mathcal{A}_{\ell}}
\underbrace{\left( \begin{array}{c}
{\bf u}_{\ell} \\
{\bf p}_{\ell}
\end{array}\right)}_{x_{\ell}} = b_{\ell},
\end{eqnarray*}
where $\mathcal{A}_{\ell}\in \mathbb{R}^{(n_{\ell}+m_{\ell}) \times (n_{\ell}+m_{\ell})}$ is a symmetric matrix, and $n_{\ell}$ and $m_{\ell}$ are the dimensions of finite-element spaces ${\bf V}_{h,\ell}$ and $Q_{h,\ell}$, respectively. We note that $A_{\ell}\in \mathbb{R}^{n_{\ell}\times n_{\ell}}$ is a symmetric and positive-definite matrix (when $\alpha$ is sufficiently large for the ${\bf H}(\text{div})$-conforming case).

The monolithic multigrid methods are defined based on the following components:
\begin{itemize} 
\item \emph{The grid transfer operations}: For the prolongation and restriction of vectors (or corresponding finite-element functions) between different levels, we use the canonical operators. The coupled prolongation
between levels ${\ell}-1$ and ${\ell}$ is given by
\begin{eqnarray*}
\mathcal{P}_{\ell} = \left(\begin{array}{cc}
P^{\bf V}_\ell & 0\\
0            & P^Q_\ell
\end{array}\right),
\end{eqnarray*}
where $P^{\bf V}_\ell\in \mathbb{R}^{n_{\ell} \times n_{{\ell}-1}}$ is the matrix representation of the finite-element interpolation operator associated with the natural embedding ${\bf V}_{h,{\ell}-1}\subset {\bf V}_{h,\ell}$ and $P^Q_\ell\in \mathbb{R}^{m_{\ell} \times m_{{\ell}-1}}$ is the matrix representation of the finite-element interpolation operator corresponding to the natural embedding $Q_{h,{\ell}-1}\subset Q_{h,\ell}$.  
For the restriction operator, $\mathcal{R}_{\ell}$, between levels ${\ell}$ and ${\ell}-1$, we take $\mathcal{R}_{\ell} = \mathcal{P}_{\ell}^T$.

We note that implementing interpolation and restriction through the natural embeddings for both continuous and discontinuous Lagrange spaces is straightforward due to their use of point evaluation dual basis functionals.  For ${\bf H}(\text{div})$-conforming spaces, integral dual basis functionals are more common~\cite{FenicsFEM}.  In the version of Firedrake~\cite{Firedrake} used in this study, grid transfers for ${\bf H}(\text{div})$-conforming spaces are implemented by first embedding the functions in appropriate discontinuous Lagrange spaces, transferred there, then projected (without loss of information) back into the appropriate space.

\item \emph{Coarse-grid operators}:  In all of the experiments in this paper, we use rediscretization to define the operator $\mathcal{A}_{\ell}$ for all levels of the multigrid hierarchy.  For the Taylor-Hood discretizations, this is equivalent to using Galerkin coarsening, since the discretization is conforming.  For the ${\bf H}(\text{div})$-conforming spaces, this is not equivalent, because of the scaling by $h_e^{-1}$ on the penalty term.  A comparison of Galerkin and rediscretization coarsening for the lowest-order ${\bf H}(\text{div})$-conforming case has been considered in the literature~\cite{ElementWiseExtendedVanka}, where rediscretization was seen to outperform Galerkin coarsening in some instances.

\item \emph{Relaxation}:  In the following section, we extend Vanka relaxation schemes from the literature to both the Taylor-Hood and ${\bf H}(\text{div})$-conforming cases, explaining in detail both the algorithm and the motivation for the underlying choice of Vanka patches.  To enable parallel computation, we consider only additive variants of Vanka relaxation~\cite{AdditiveVanka,PFarrell_etal_2019a}.  Within the multigrid V-cycle, we always accelerate these additive schemes using the Chebyshev iteration.  To define the Chebyshev polynomials, we let $\lambda$ be an estimate of the largest eigenvalue of the relaxation-preconditioned linear system on any level, obtained from the standard GMRES estimation. Then, we define the lower and upper bounds of this interval to be $\frac{1}{4}\lambda$ and $1.1\lambda$, respectively.  The lower bound is determined because we wish to relax the upper $3/4$ of the spectrum on each level (leaving the lowest $1/4$ of the spectrum to be attenuated by coarse-grid correction).  The upper bound allows a ``safety factor'' in case $\lambda$ is an underestimate of the largest eigenvalue.
\end{itemize}
%%%%%%%%%%%%%%%%%%%%%%%%%%%%%%%%%
\subsection{Space decomposition and patch-based relaxation}
%%%%%%%%%%%%%%%%%%%%%%%%%%%%%%%%%

In what follows, we drop the subscript corresponding to the level within the multigrid hierarchy, since the resulting algorithms are the same on every level.
Given a product space ${\bf V}_h\times Q_h$, a space decomposition of ${\bf V}_h\times Q_h$ takes the form of 
\begin{eqnarray}\label{eq:spacedecomp}
{\bf V}_h\times Q_h = \sum\limits_{i=1}^{J} \left({\bf V}_h^{(i)}\times Q_h^{(i)}\right),
\end{eqnarray}
meaning that every $({\bf v}_h,q_h)\in {\bf V}_h\times Q_h$ can be written as $({\bf v}_h,q_h) = \sum\limits_{i = 1}^{J}({\bf v}_h^{(i)},q_h^{(i)})$, for ${\bf v}_h^{(i)}\in {\bf V}_h^{(i)}$ and $q_h^{(i)}\in Q_h^{(i)}$. This decomposition is not necessarily unique.  Theory and practice of solvers and relaxation schemes based on such decompositions are common in the literature~\cite{Spacedecomp, PCPATCH}.  At their core, subspace decomposition solvers solve a single problem for each subspace in~\eqref{eq:spacedecomp}, and they differ in whether they use an additive (or parallel) approach, or a multiplicative (sequential) approach, in which the residual is updated after each subdomain is visited.   In this paper, we only consider the additive approach, which is presented in Algorithm~\ref{alg:ParSubCorrect}.  We note that the differences between applying Algorithm~\ref{alg:ParSubCorrect} to the Taylor-Hood and ${\bf H}(\text{div})$-conforming cases are in the choice of space decomposition, detailed below, and in the use of the conforming bilinear form $a(\cdot,\cdot)$ from~\eqref{eq:VariationalStokesa} versus the non-conforming $a_h(\cdot,\cdot)$ from~\eqref{eq:DiscontBilinearForm}.

%%%%%%%%%%%%%%%%%%%%%%%%%%%%%%%%%
\begin{algorithm}[t]%[]
\caption{ {\bf Additive (parallel) subspace correction method}}
\label{alg:ParSubCorrect}
%\algsetup{linenodelimiter=.}
%\begin{footnotesize}
{\bf Inputs}: 
\begin{itemize}
\item initial approximation $(\hat{\bf u}_h, \hat{p}_h)\in {\bf V}_h\times Q_h$
\item the space decomposition ${\bf V}_h\times Q_h = \sum\limits_{i=1}^{J} \left({\bf V}_h^{(i)}\times Q_h^{(i)}\right)$
\item the weighting operators $w_i: {\bf V}_h^{(i)}\rightarrow {\bf V}_h^{(i)}$ and ${\gamma}_i: Q_h^{(i)}\rightarrow Q_h^{(i)}$
\item bilinear forms $a(\cdot,\cdot)$ and $b(\cdot,\cdot)$ and the external body force ${\bf f}:\Omega\rightarrow \mathbb{R}^2$
\end{itemize}
{\bf Output}: Updated approximation $(\hat{\bf u}_h, \hat{p}_h)\in {\bf V}_h\times Q_h$
\begin{algorithmic}[1]
\FOR {$i=1$ to $J$}
\STATE Find the solution  $\left(\delta {\bf u}^{(i)},\delta p^{(i)}\right)\in {\bf V}_h^{(i)}\times Q_h^{(i)}$ of the local problem 
\begin{eqnarray*}
\left\{\begin{array}{l}
a(\delta {\bf u}^{(i)},{\bf v}_h) + b({\bf v}_h,\delta p^{(i)}) = \langle{\bf f},{\bf v}_h\rangle - a(\hat{\bf u}_h,{\bf v}_h) - b({\bf v}_h,\hat{p}_h), \\
           b(\delta {\bf u}^{(i)},q_h) = - b(\hat{\bf u}_h,q_h),
\end{array}
\right.
\end{eqnarray*}
\hspace{12pt} for all $({\bf v}_h, q_h) \in {\bf V}_h^{(i)}\times Q_h^{(i)}$.
\ENDFOR
\STATE $(\hat{\bf u}_h, \hat{p}_h) \leftarrow \left(\hat{\bf u}_h + \sum\limits_{i=1}^{J} w_i(\delta {\bf u}^{(i)}), \hat{p}_h + \sum\limits_{i=1}^{J} {\gamma}_i(\delta p^{(i)})\right)$ 
\STATE \Return $(\hat{\bf u}_h,\hat{p}_h)$
\end{algorithmic}
%\end{footnotesize}
\end{algorithm}
%%%%%%%%%%%%%%%%%%%%%%%%%%%%%%%%%

We are interested in the use of subspace correction methods as relaxation schemes in multigrid solvers.  The use of such relaxation schemes dates back to the pioneering work of Vanka~\cite{VankaType}, who considered the marker-and-cell (MAC) scheme finite-difference discretization for the Navier-Stokes equations, defining patches element-wise for a discretization with pressures located at element centers and normal velocity DoFs on element edges.  This was extended to the Crouzeix-Raviart finite-element discretization by John and Tobiska~\cite{JohnVankaTypeNavierStokes}, again using element-wise patches.  The first multigrid schemes for a discretization without element-centered pressure DoFs were proposed by John and Matthies~\cite{JohnPatchesStokes}.  In that work, they considered both element-wise patches and ``pressure node oriented'' relaxation, writing (in the notation above) that one should select the patches so that ``dim$(Q_h^{(i)}) = 1$ and ${\bf V}_h^{(i)}$ is the set of velocity degrees of freedom which are connected to the pressure degree of freedom of $Q_h^{(i)}$ in the matrix $B$''.  This success was mimicked in the so-called ``algebraic Vanka'' approaches~\cite{JAdler_etal_2014c,TBenson_etal_2015a}.  More recently, a topological patch construction approach was implemented in PCPATCH~\cite{PCPATCH}.  In this work, we explore a topological variant of John and Matthies' ``pressure node oriented'' relaxation, as explained below.

Firedrake~\cite{Firedrake} makes use of the DMPlex~\cite{doi:10.1137/15M1026092} data management structure in PETSc~\cite{petsc-web-page, petsc-user-ref}.  DMPlex provides standard data structures and operations for working with unstructured grids, associating discrete DoFs with nodes, edges, faces, and cells of the grid, and providing operations for finding adjacency relations between topological objects.  Thus, in Firedrake, we can naturally specify patches by their topological construction.
For example, to construct the set of vertex-star patches, we loop over all of the vertices of the mesh and apply the star operation~\cite[Chapter 1, Section 2]{CellPatch} to gather a set of entities adjacent to each vertex. The star of a vertex is defined by the (interiors of all) edges, faces, and elements incident on the vertex, along with the vertex itself.  PETSc's DMPlex data structure provides straightforward tools to construct the topological star, and then to gather the discretized degrees of freedom on each topological object contained within the star.

\begin{figure}[t]
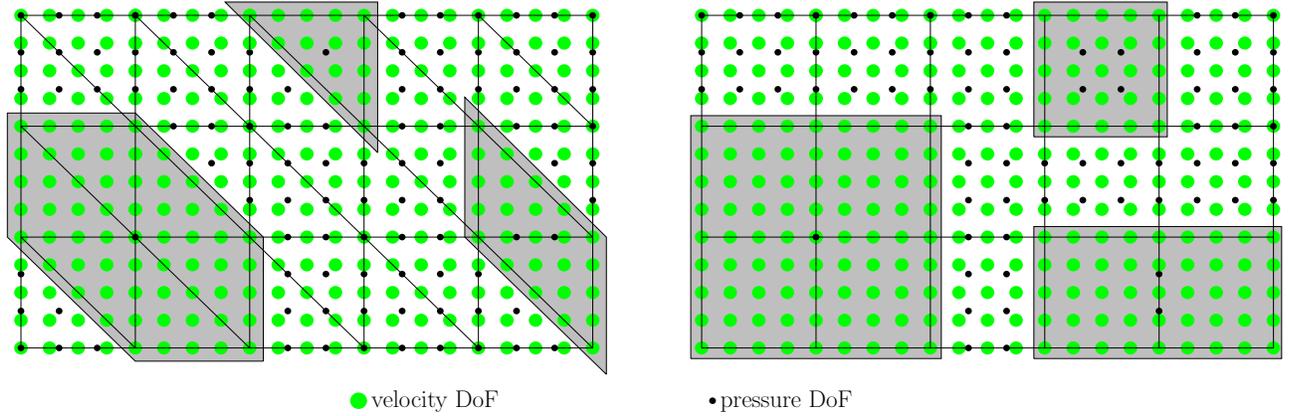
 
\begin{center}

\end{center}
\caption{{\bf Left}: Topological Vanka patches for $({\bf P_4},P_3)$ Taylor-Hood discretization, showing nodal (left), edge (right) and cell (top) patches for composite Vanka relaxation scheme.
  {\bf Right}: Topological Vanka patches for $({\bf Q_4},Q_3)$ Taylor-Hood discretization, showing nodal (left), edge (right) and cell (top) patches for composite Vanka relaxation scheme. Big (green) and small (black) circles are associated with velocity and pressure DoFs, respectively.  Within each patch, we show only the included degrees of freedom, suppressing pressure DoFs that are not included in the patches.}
\label{fig:CompositeVankaPatchStokes}
\end{figure}

For Taylor-Hood elements, the pressure degrees of freedom live on a variety of topological entities in two dimensions, depending on the order of the discretization.  Classical Vanka relaxation for the $({\bf P_2},P_1)$ and $({\bf Q_2},Q_1)$ discretizations involves only the vertex-based pressures, so we can form suitable patches by taking the pressure degree of freedom at each node, and the velocity degrees of freedom on the closure of the associated vertex-star patch.  Depending on symmetries in the mesh (that may lead to ``coincidental'' zero entries in the stencil of the matrix), this generally matches the ``pressure node oriented'' Vanka of John and Matthies.  Similarly, for $({\bf P_3},P_2)$ and $({\bf Q_3},Q_2)$, we can form patches that match the ``pressure node oriented'' patches up to coincidental zero entries.  For $({\bf P_3},P_2)$, there are both nodal and edge-based pressure DoFs and, so, we construct two types of patches, one that matches the construction of the vertex-based Vanka patches (with the pressure DoF at each vertex and all velocity DoFs in the closure of the vertex star), and one that is formed by taking each edge-based pressure DoF and all velocity DoFs in the closure of the edge star, meaning those velocity DoFs on the edge or its incident nodes, and in the closure of the two elements adjacent to the edge.  For $({\bf Q_3},Q_2)$, we take three types of patches, using vertex and edge Vanka patches, as for $({\bf P_3},P_2)$, but adding patches that are cell-based, since $Q_2$ has a cell-centered pressure DoF.  These cell-based Vanka patches consist of the pressure DoF located within a cell as well as all velocity DoFs on the closure of the cell (in 2D, the cell itself, and all DoFs on the edges and vertices of the cell).  For $k\geq 4$, we start to have topological entities with more than one pressure DoF in $P_{k-1}$, since there are $k-2$ pressure DoFs on each edge, and $(k-2)(k-3)/2$ pressure DoFs in the cell.  Thus, for $k\geq 4$, we depart from the ``pressure node oriented'' patches to consider ``topological Vanka'' patches that contain 1 or more pressure DoFs on the same topological entity, defining $Q_h^{(i)}$, and all velocity DoFs on the closure of the star of that entity.  Since our construction makes use of multiple patches of different type, we refer to this as a composite Vanka relaxation scheme.  Sample patches for $k=4$ are shown in Figure~\ref{fig:CompositeVankaPatchStokes}.  We note that the triangular case gives patches with consistently fewer degrees of freedom than the quadrilateral case.

%%%%%%%%%%%%%%%%%%%%%%%%%%%%%%%

\begin{figure}[t]
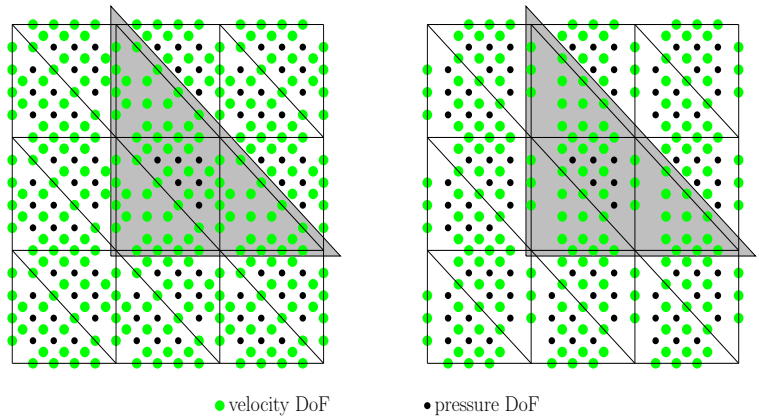
 
\begin{center}

\end{center}
\caption{{\bf Left:} Extended Vanka patches for $(\mathcal{ BDM}_3,{dP}_2)$.  {\bf Right:} Extended Vanka patches for $(\mathcal{RT}_3,{dP}_2)$.  Big (green) and small (black) circles are associated with velocity and pressure DoFs, respectively. Within each patch, we show only the included degrees of freedom, suppressing pressure DoFs that are not included in the patches.}
\label{fig:ExtendVankaPatchBDMDGStokes}
\end{figure}

The use of discontinuous pressure spaces for the ${\bf H}(\text{div})$-$L^2$ formulation greatly simplifies the choice of Vanka patches, since there are only cell-based pressure degrees of freedom.  Following the original prescriptions for patch selection~\cite{VankaType, JohnVankaTypeNavierStokes, JohnPatchesStokes}, we are tempted to choose to use only the cell-based patches above for these cases.  However, Adler et al.~\cite{ElementWiseExtendedVanka} demonstrate that, for the lowest-order $(\mathcal{BDM}_1,dP_0)$ discretization, these patches do not lead to the most effective solvers.  Instead, here we adopt the same extended cell-based Vanka patches proposed there.  For each triangle in the mesh, we construct Vanka patches so that all pressure DoFs on the triangle are included in the patch, along with all velocity DoFs on the element (the cell, plus its adjacent edges and vertices) and on the closure of the three elements that share edges with the cell.
Figure~\ref{fig:ExtendVankaPatchBDMDGStokes} shows these patches for both cases, with $(\mathcal{BDM}_3,dP_2)$ elements and $(\mathcal{RT}_3,dP_2)$ elements.

%%%%%%%%%%%%%%%%%%%%%%%%%%%%%%%
%%%%%%%%%%%%%%%%%%%%%%%%%%%%%%%%%%%%%

Firedrake provides two implementations of patch-based relaxation schemes.  PCPATCH~\cite{PCPATCH} provides a callback-based implementation, where the global solve can be done in a matrix-free way.  For each patch (on each multigrid level), the elements needed to assemble the local operator are identified and passed to Firedrake, which assembles and stores the patch matrices (or their inverses or factored forms) independently of the global operator.  ASMPatchPC, in contrast, is based on an assembled global operator (i.e., a matrix-full implementation).  In this implementation, we identify the global indices of each DoF that appears in each patch, and extract the patch submatrix from the globally assembled matrix directly.  While PCPATCH has some advantages from not requiring the global matrix to be assembled, there are also disadvantages to this approach.  One key disadvantage is that PCPATCH requires redundant computation when assembling overlapping patches, since each patch subproblem is assembled independently.  A second disadvantage is that PCPATCH is harder to extend to new patch constructions than ASMPatchPC; this is particularly relevant here, since the extended cell Vanka patches that we construct for ${\bf H}(\text{div})$-$L^2$ discretizations are not currently implemented in PCPATCH.  Combined with this is the fact that there is no current support for facet integrals in assembling the local problems in PCPATCH, so we cannot apply it to discretizations involving the weak form in~\eqref{eq:DiscontBilinearForm}.  In the numerical results below, we compare the efficiency of PCPATCH and ASMPatchPC for the Taylor-Hood discretizations, but note that no such comparison is possible for the other discretizations considered here.

%%%%%%%%%%%%%%%%%%%%%%%
\section{Numerical Experiments}\label{sec:numerics}
All numerical results in this paper are based on a manufactured solution to the Stokes problem posed on $\Omega = [0,1]^2$, given by
\begin{equation}\label{eq:model_problem}
{\bf u}(x,y) = \left( \begin{array}{c}
\sin(\pi x)\cos(\pi y) \\
-\cos(\pi x)\sin(\pi y)
\end{array}\right) \text{ and }  p(x,y) = 0,
\end{equation}
with boundary conditions and right-hand side function ${\bf f}$ chosen give this analytical solution.  All numerical results were computed on a workstation with two 8-core 1.7GHz Intel Xeon Bronze 3106 CPUs and 384 GB of RAM.  All timings are reported in seconds for runs on 8 cores.

We present results for this problem in three parts, considering the Taylor-Hood discretization on triangles in Section~\ref{ssec:TH_triangle}, the Taylor-Hood discretization on quadrilaterals in Section~\ref{ssec:TH_quad}, and the ${\bf H}(\text{div})$-$L^2$ discretizations on triangles in Section~\ref{ssec:Hdiv_DG_triangle}.  We use $\mathcal{TH}$ to denote results for the Taylor-Hood discretizations and $\mathcal{BDM}$ or $\mathcal{RT}$ to denote results for the ${\bf H}(\text{div})$-$L^2$ discretizations.  In all cases, we use a subscript to denote the type of grid considered, with $\mathcal{TH}_T$ denoting a Taylor-Hood discretization on a triangular grid, and $\mathcal{TH}_Q$ denoting a Taylor-Hood discretization on a quadrilateral grid.  We also use superscripts to differentiate between solvers using PCPATCH, denoted as $\mathcal{TH}_T^{PC}$, and those using ASMPatchPC, denoted as $\mathcal{TH}_T^{AS}$.  Since the ${\bf H}(\text{div})$-$L^2$ discretizations are only considered on triangles and using ASMPatchPC, we note this notation is superfluous, but include it for consistency.  For each method, we will consider variation with respect to three parameters, the order of the discretization, $k$ (which we take to be the order of the velocity space), the number of pre- and post-relaxation sweeps on each level, $\nu$, and the number of levels in the hierarchy, $\ell$.  In all experiments, we use a $5\times 5$ mesh of the unit square as the coarsest grid, with square elements cut from top left to bottom right to form two triangular elements when using triangular grids.  Thus, we generally refer to a method by $\mathcal{TH}_T^{PC}(k,\nu,\ell)$, specifying one or more of the parameters as needed.  The initial guess for the iterative solver is taken to be the zero vector, and the solver is deemed to have converged when the discrete Euclidean norm of the residual has been reduced by a relative factor of $10^{10}$.  We consider unrestarted FGMRES as the outer solver, allowing a maximum of 100 iterations, unless otherwise noted.  We more closely examine stopping criteria and discretization error in Section~\ref{ssec:stopping}.

%Fig4%%%%%%%%%%%%%%%%%%%%%%%%%%%%
%%%%%%%%%%%%%%%%%%%%%%%%%%%%%%%
\begin{figure}[t]
\graphicspath{{./Figures/}}
\begin{center}
\includegraphics{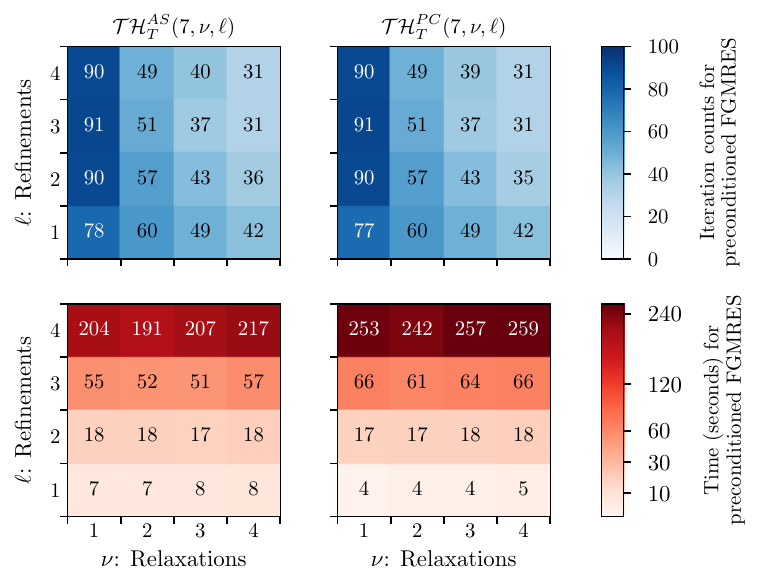}
\end{center}
\caption{Comparison of iteration counts ({\bf top row}) and time-to-solution ({\bf bottom row}) for the Taylor-Hood solvers on triangular grids with order $k=7$.  Results at {\bf left} use ASMPatchPC for the Vanka relaxation, while those at {\bf right} use PCPATCH.}
\label{fig:P7P6ASMPatchPCvsPCPATCH_It_Time_FGMRES}
\end{figure}
%%%%%%%%%%%%%%%%%%%%%%%%%%%%%%%
%%%%%%%%%%%%%%%%%%%%%%%%%%%%%%%

%%%%%%%%%%%%%%%%%%%%%%%
%%%%%%%%%%%%%%%%%%%%%%%
\subsection{Taylor-Hood on triangular grids}\label{ssec:TH_triangle}

We first consider the performance of the monolithic multigrid method for fixed polynomial order $k=7$ with $1\leq \nu \leq 4$ and $1 \leq \ell \leq 4$.  Figure~\ref{fig:P7P6ASMPatchPCvsPCPATCH_It_Time_FGMRES} presents iteration counts to convergence (top) and time-to-solution (bottom) for the solvers implemented using ASMPatchPC (left) and PCPATCH (right).  As expected, the iteration counts for the two implementations are essentially identical, with slight differences that can be attributed to small differences in floating-point calculations.  Further, we see that the iterations counts for $V(1,1)$ cycles (when $\nu = 1$) are quite high, with a slight increase as we go from a two-level scheme with $\ell=1$ refinement of the base mesh to multilevel solvers on finer grids.  In contrast, when $\nu > 1$, we see somewhat more reasonable iteration counts, generally decreasing as we go to finer grids.

Considering the CPU times reported in the bottom part of Figure~\ref{fig:P7P6ASMPatchPCvsPCPATCH_It_Time_FGMRES}, we notice two things.  First of all, the solvers implemented using ASMPatchPC are generally about 20\% faster on finer grids, while those using PCPATCH are somewhat faster only when $\ell = 1$.  Secondly, in both cases, the CPU times scale consistently with the total number of degrees of freedom in the system, growing roughly by a factor of 4 when the finest grid is refined once.  What is quite notable in these results is that even though there is a lot of variation in the iteration counts as we vary the number of relaxations, $\nu$, there is relatively little variation in the total time-to-solution for either implementation (even though relaxation takes the dominant amount of time in these simulations).  This is almost certainly because the reduction in iterations as we increase $\nu$ offsets the added cost per iteration.  Nonetheless, in both cases, we see some advantage to selecting $\nu=2$, which seems to balance the cost per iteration versus the number of iterations to convergence on finer grids, saving about 5\% of the CPU time.

Figure~\ref{fig:All_Timings_PP_ASMPatchPCvsPCPATCH} presents CPU time-to-solution for the two variants (using ASMPatchPC at left and PCPATCH at right) with $\ell= 4$, as we vary both order, $k$, and number of pre- and post-relaxation sweeps per cycle, $\nu$.  Here, we see generally comparable results to those shown for $k=7$ in Figure~\ref{fig:P7P6ASMPatchPCvsPCPATCH_It_Time_FGMRES}.  The ASMPatchPC implementation is generally about 20\% faster than the PCPATCH implementation, and choosing $\nu = 2$ generally yields the fastest solver for a fixed value of $k$, although the overall benefit can be small.  Comparing time-to-solution for $k=4$ and $k=8$, we see an increase in time-to-solution by a factor of roughly 8, suggesting a scaling of $\mathcal{O}(k^3)$ for time-to-solution.  Since (in 2D), the number of DoFs in each patch system scales roughly like $\mathcal{O}(k^2)$, this strongly outperforms a pessimistic estimate of the cost of (dense) patch solves, which would be $\mathcal{O}(k^6)$.  Some of this improvement may simply be due to the relatively small orders considered here (noting that a smaller factor of increase is seen going from $k=3$ to $k=6$), and some may be due to advantages due to structural or coincidental zero entries that appear in the patch systems.  We note, in fact, that this scaling also outperforms our expectation for that of simply forming the matrix or computing a matrix-vector product, both of which should scale like $\mathcal{O}(k^4)$, since each element stiffness matrix contains dense regions over $\mathcal{O}(k^2)$ DoFs.

Finally, in Figure~\ref{fig:All_PP_ASMPatchPC_ItFGMRES}, we consider the robustness of the solvers to polynomial order, focusing on the case of using ASMPatchPC as the relaxation implementation.  Here, we first notice that the solver is far from scalable with small values of $\nu$, showing substantial increases in the number of iterations to convergence (and even failing to converge) as we increase $k$ from 3 to 8 with $\nu = 1$, even for coarse grids with levels in the multigrid hierarchy.  In contrast, for the finest grids with $\nu = 3$ or $4$, the increase in iterations is less pronounced, with the iteration counts for $\nu = 4$ and $\ell = 4$ increasing only from $26$ for $k=3$ to $35$ for $k=8$, \revise{showing, at best, a weak sort of $p$-robustness}.  As seen above, however, even though the iteration counts are more scalable with large $\nu$, the CPU times still are minimized by consistently taking $\nu = 2$, where there is about an 80\% growth in iteration counts at $\ell = 4$ when moving from $k=3$ to $k=8$.  Moreover, we note that for any fixed polynomial order, $k$, if $\nu > 1$, we see decreasing iteration counts as the number of levels in the multigrid hierarchy increases, as we would hope to see from \revise{an $h$-robust} multigrid solver.
\subsection{Taylor-Hood on quadrilateral grids}\label{ssec:TH_quad}

\begin{figure}[t]
\graphicspath{{./Figures/}}
\begin{center}
\includegraphics{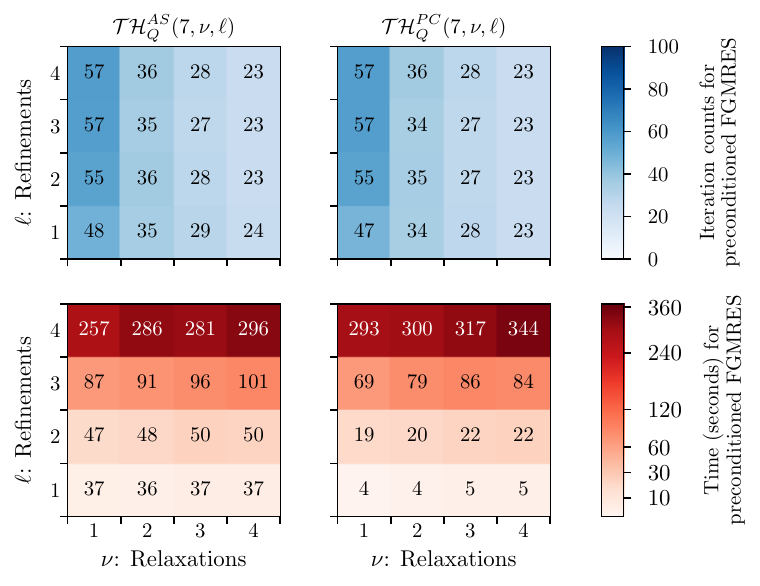}
\end{center}
\caption{Comparison of iteration counts ({\bf top row}) and time-to-solution ({\bf bottom row}) for the Taylor-Hood solvers on quadrilateral grids with order $k=7$.  Results at {\bf left} use ASMPatchPC for the Vanka relaxation, while those at {\bf right} use PCPATCH.}
\label{fig:Q7Q6_ASMPatchPCvsPCPATCH_It_Time_FGMRES}
\end{figure}

Even though we consider only regular grids, it is common to see substantial differences in performance between solvers applied to equivalent discretizations on triangular or quadrilateral meshes.  Thus, we now repeat the experiments from the previous section but with quadrilateral meshes.  Figure~\ref{fig:Q7Q6_ASMPatchPCvsPCPATCH_It_Time_FGMRES} performs the same comparison as in Figure~ \ref{fig:P7P6ASMPatchPCvsPCPATCH_It_Time_FGMRES}, but for the case of Taylor-Hood discretizations on quadrilateral grids.  As before and as expected, we see only trivial differences in iteration counts between the implementations making use of ASMPatchPC and PCPATCH.  Comparing CPU times between the two implementations, we see substantial speedup using PCPATCH for the cases of $\ell = 1$ and $2$ but, by the time we reach $\ell=4$, we see that ASMPatchPC again is consistently faster than PCPATCH, albeit by a smaller margin.
%Fig5%%%%%%%%%%%%%%%%%%%%%%%%%%%%
%%%%%%%%%%%%%%%%%%%%%%%%%%%%%%%
\begin{figure}[t]
\graphicspath{{./Figures/}}
\begin{center}
\includegraphics{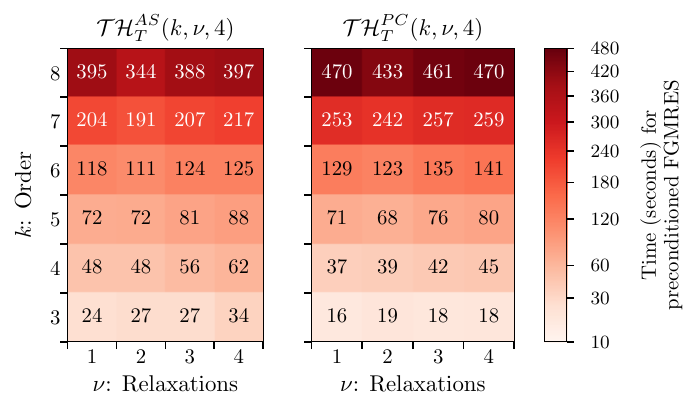}
\end{center}
\caption{Comparison of time-to-solution for two implementations of solvers for the Taylor-Hood discretization on triangular grids.  At {\bf left}, we use ASMPatchPC for relaxation.  At {\bf right}, we use PCPATCH for relaxation.}
\label{fig:All_Timings_PP_ASMPatchPCvsPCPATCH}
\end{figure}
%%%%%%%%%%%%%%%%%%%%%%%%%%%%%%%
%%%%%%%%%%%%%%%%%%%%%%%%%%%%%%%

%Fig6%%%%%%%%%%%%%%%%%%%%%%%%%%%%
%%%%%%%%%%%%%%%%%%%%%%%%%%%%%%%
\begin{figure}[H]
\graphicspath{{./Figures/}}
\begin{center}
\includegraphics{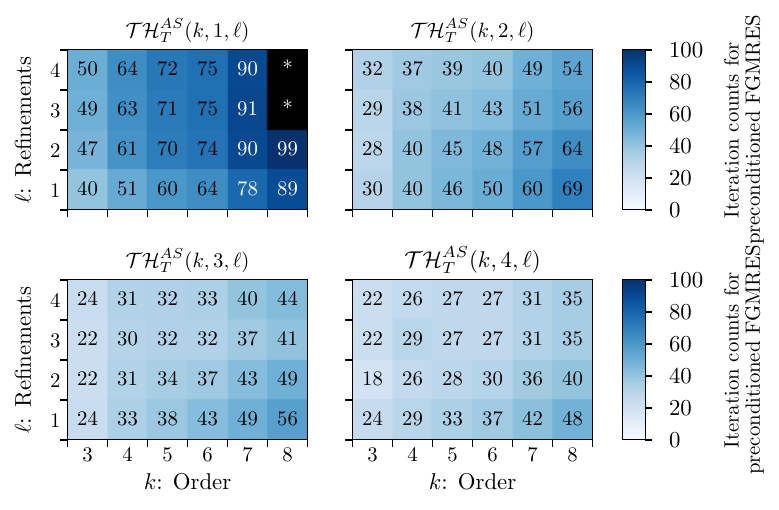}
\end{center}
\caption{Iteration counts for solving the Taylor-Hood discretization on triangular grids with varying order $k$ using ASMPatchPC for the relaxation implementation. Results indicated by a black square with white asterisk indicate that the solver did not converge to the desired tolerance within 100 iterations.}
\label{fig:All_PP_ASMPatchPC_ItFGMRES}
\end{figure}
%%%%%%%%%%%%%%%%%%%%%%%%%%%%%%%
%%%%%%%%%%%%%%%%%%%%%%%%%%%%%%%

%Fig7%%%%%%%%%%%%%%%%%%%%%%%%%%%%
%%%%%%%%%%%%%%%%%%%%%%%%%%%%%%%
%%%%%%%%%%%%%%%%%%%%%%%%%%%%%%%
%%%%%%%%%%%%%%%%%%%%%%%%%%%%%%%

\begin{figure}[t]
\graphicspath{{./Figures/}}
\begin{center}
\includegraphics{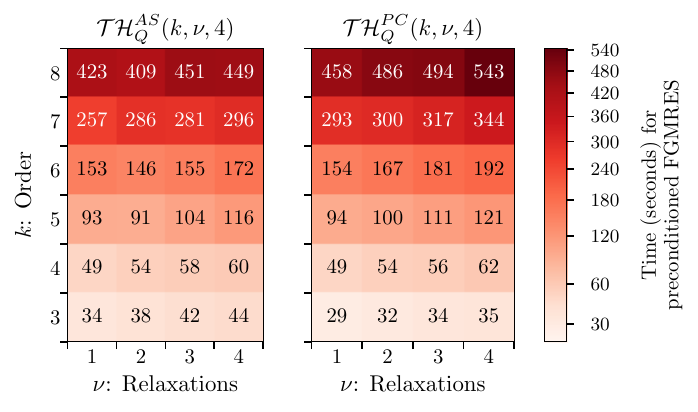}
\end{center}
\caption{Comparison of time-to-solution for two implementations of solvers for the Taylor-Hood discretization on quadrilateral grids.  At {\bf left}, we use ASMPatchPC for relaxation.  At {\bf right}, we use PCPATCH for relaxation.}
\label{fig:All_Timings_QQ_ASMPatchPCvsPCPATCH}
\end{figure}

The much more substantial difference seen in Figure~\ref{fig:Q7Q6_ASMPatchPCvsPCPATCH_It_Time_FGMRES} is in the iteration counts in the first place, with much more scalable counts shown for the case of $\nu = 1$, and more modest iteration counts for larger values of $\nu$.  Put simply, the composite Vanka solver appears to perform better on quadrilateral meshes than it does on triangular ones.  One important consequence of this is that the solver performance is optimized by taking $\nu = 1$ consistently across all refinement levels, $\ell$.  A similar story is told in Figure~\ref{fig:All_Timings_QQ_ASMPatchPCvsPCPATCH}, which shows CPU timings for the two implementations as we vary $k$ for $\ell = 4$.  Here, the most efficient solver is almost always when $\nu = 1$, with occasional slight improvements when $\nu = 2$.  For the case of $k=8$, we see about a 10\% improvement in the best time-to-solution when using ASMPatchPC than when using PCPATCH.  Considering the CPU time scaling with polynomial order, $k$, we see very similar behaviour as for the triangular grid case, with time-to-solution growing by about a factor of 8 when moving from $k=4$ to $k=8$, suggesting scaling like $\mathcal{O}(k^3)$, which is substantially better than expected (and, indeed, comparable to that of sum factorization algorithms~\cite{SumFactorization}).  A natural consideration for future work in this case is to see if the recently proposed fast diagonalization preconditioners of Brubeck and Farrell~\cite{FDM} can be extended to give comparable convergence but with improved cost guarantees.

\begin{figure}[H]
\graphicspath{{./Figures/}}
\begin{center}
\includegraphics{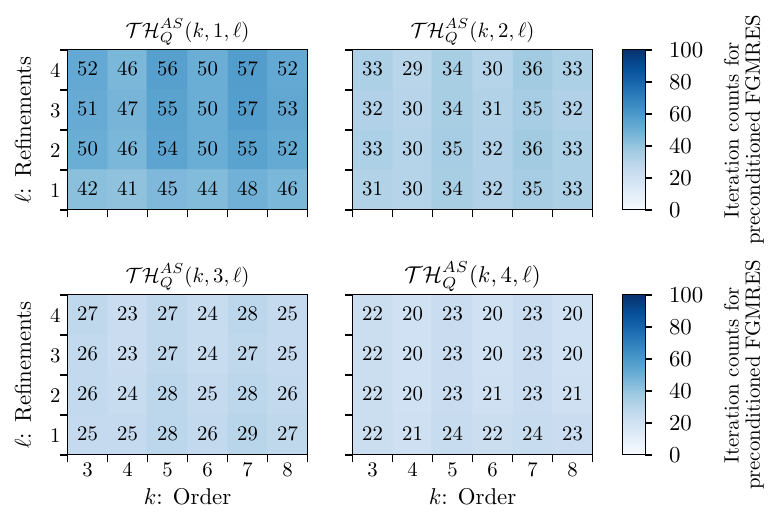}
\end{center}
\caption{Iteration counts for solving the Taylor-Hood discretization on quadrilateral grids with varying order $k$ using ASMPatchPC for the relaxation implementation.}
\label{fig:All_QQ_ASMPatchPC_ItFGMRES}
\end{figure}

Finally, Figure~\ref{fig:All_QQ_ASMPatchPC_ItFGMRES} shows iteration counts for the solver for Taylor-Hood discretizations on quadrilateral grids using ASMPatchPC for the relaxation.  Here, in contrast to Figure~\ref{fig:All_PP_ASMPatchPC_ItFGMRES}, we see almost perfectly scalable results with respect to both grid refinement and polynomial order, \revise{showing both $h$- and $p$-robustness of the solver (and, consequently, $hp$-robustness)}.  While there remains some variation in these iteration counts, we see that they exhibit almost none of the consistent growth seen for triangular grids with increasing polynomial order.  Thus, while performance with $k=3$ is comparable between the two types of grids, we see notably improved performance for quadrilateral grids with larger values of $k$.  This is particularly interesting when comparing the timings shown for the two cases in Figures~\ref{fig:All_Timings_PP_ASMPatchPCvsPCPATCH} and~\ref{fig:All_Timings_QQ_ASMPatchPCvsPCPATCH}, which show that the solvers on triangular grids are still cheaper than those on quadrilateral grids, despite this growth in iteration counts.  We suspect this comes from the larger sizes of the patches for quadrilateral grids compared to the triangular case, so that each composite Vanka relaxation sweep is more expensive, leading to a higher cost per iteration.
%Fig8%%%%%%%%%%%%%%%%%%%%%%%%%%%%
%%%%%%%%%%%%%%%%%%%%%%%%%%%%%%%
%%%%%%%%%%%%%%%%%%%%%%%%%%%%%%%
%%%%%%%%%%%%%%%%%%%%%%%%%%%%%%%

%Fig9%%%%%%%%%%%%%%%%%%%%%%%%%%%%
%%%%%%%%%%%%%%%%%%%%%%%%%%%%%%%
%%%%%%%%%%%%%%%%%%%%%%%%%%%%%%%
%%%%%%%%%%%%%%%%%%%%%%%%%%%%%%%

\subsection{${\bf H}(\text{div})$-$L^2$ discretizations on triangular grids}\label{ssec:Hdiv_DG_triangle}

\begin{figure}[t]
\graphicspath{{./Figures/}}
\begin{center}
\includegraphics{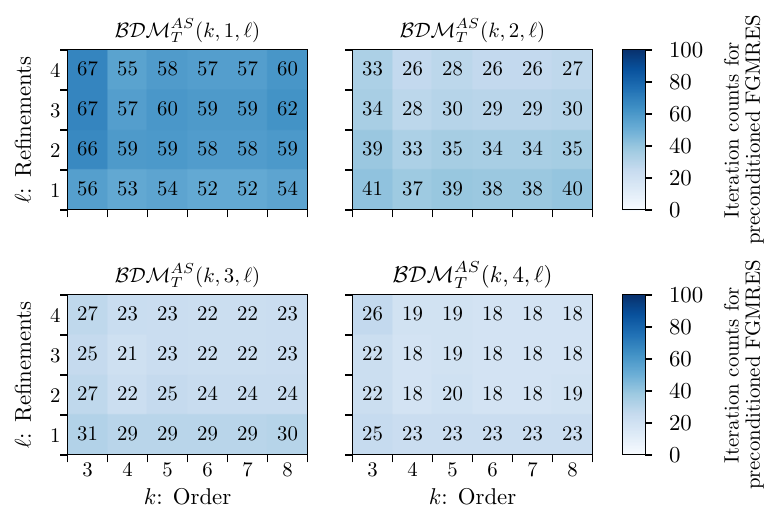}
\end{center}
\caption{Iteration counts for solving the $(\mathcal{BDM}_k,dP_{k-1})$ discretization on triangular grids with varying order $k$ using ASMPatchPC for the relaxation implementation.}
\label{fig:All_BDMdP_ASMPatchPC_ItFGMRES}
\end{figure}

Figure~\ref{fig:All_BDMdP_ASMPatchPC_ItFGMRES} presents results corresponding to Figures~\ref{fig:All_PP_ASMPatchPC_ItFGMRES} and~\ref{fig:All_QQ_ASMPatchPC_ItFGMRES}, but for the $(\mathcal{BDM}_k,dP_{k-1})$ discretization with varying orders, $k$.  Here, we use only ASMPatchPC to implement the relaxation, because the extended Vanka patches shown in Figure~\ref{fig:ExtendVankaPatchBDMDGStokes} are not currently implemented in PCPATCH.  Most notable in these results is that the iteration counts appear perfectly scalable, as they were for the quadrilateral Taylor-Hood case in Figure~\ref{fig:All_QQ_ASMPatchPC_ItFGMRES}, \revise{again showing $h$-, $p$-, and $hp$-robustness}.  \revise{We note that similar experiments with cell-based patches showed a strong lack of robustness, with iteration counts growing dramatically with mesh size.}
%%%%%%%%%%%%%%%%%%%%%%%%%%%%%%%
%%%%%%%%%%%%%%%%%%%%%%%%%%%%%%%
%Fig10%%%%%%%%%%%%%%%%%%%%%%%%%%%%
%%%%%%%%%%%%%%%%%%%%%%%%%%%%%%%
%%%%%%%%%%%%%%%%%%%%%%%%%%%%%%%
%%%%%%%%%%%%%%%%%%%%%%%%%%%%%%%

\begin{figure}[t]
\graphicspath{{./Figures/}}
\begin{center}
\includegraphics{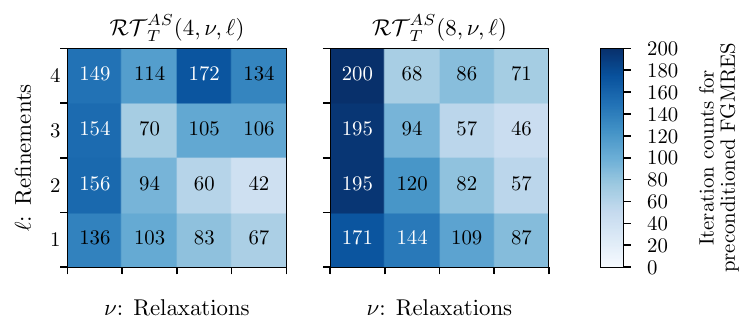}
\end{center}
\caption{Iteration counts for solving the $(\mathcal{RT}_k,dP_{k-1})$ discretization on triangular grids with $k=4$ ({\bf left}) and $k=8$ ({\bf right}) using ASMPatchPC for the relaxation implementation.  Here, FGMRES(200) is used as the Krylov method.}
\label{fig:RT4dP3_RT8dP7_ASMPatchPC_ItFGMRES}
\end{figure}

In contrast to the scalable convergence seen with the $(\mathcal{BDM}_k,dP_{k-1})$ discretization, Figure~\ref{fig:RT4dP3_RT8dP7_ASMPatchPC_ItFGMRES} shows convergence for the $(\mathcal{RT}_k,dP_{k-1})$ discretization for $k=4$ and $8$.  Here, we see very poor scalability, both with respect to increasing number of levels in the hierarchy, $\ell$, and polynomial order, $k$.  These results are, notably, much worse than those reported above for any of the other solvers considered and, consequently, we focus our experiments in this section on the $(\mathcal{BDM}_k,dP_{k-1})$ discretization.  We do note that in timing results (not reported here), the best time-to-solution for the $(\mathcal{RT}_8,dP_{7})$ discretization was almost 700 seconds, in contrast to about 500 seconds reported in the results below for the $(\mathcal{BDM}_8,dP_{7})$ discretization.  Thus, it takes almost twice as long to achieve lower-accuracy results using Raviart-Thomas elements in this setting than when using Brezzi-Douglas-Marini elements.
%Fig11%%%%%%%%%%%%%%%%%%%%%%%%%%%%
%%%%%%%%%%%%%%%%%%%%%%%%%%%%%%%
%%%%%%%%%%%%%%%%%%%%%%%%%%%%%%%
%%%%%%%%%%%%%%%%%%%%%%%%%%%%%%%

\begin{figure}[t]
\graphicspath{{./Figures/}}
\begin{center}
\includegraphics{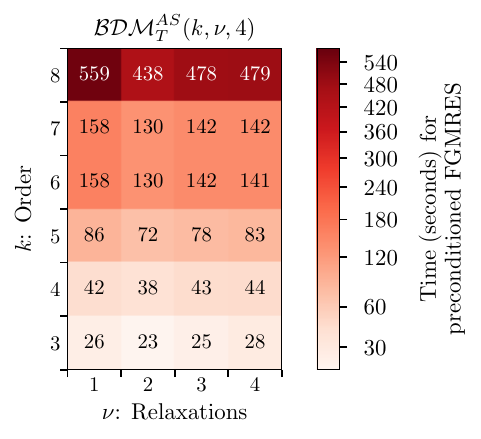}
\end{center}
\caption{Time-to-solution for the $(\mathcal{BDM}_k,dP_{k-1})$ discretization on triangular grids with varying order $k$ using ASMPatchPC for the relaxation implementation for $\ell = 4$.}
\label{fig:All_Timings_BDMdP_ASMPatchPC}
\end{figure}

Finally, in Figure~\ref{fig:All_Timings_BDMdP_ASMPatchPC}, we provide timing results for the $(\mathcal{BDM}_k,dP_{k-1})$ discretization on triangular grids with $\ell = 4$, for comparison with timings in Figures~\ref{fig:All_Timings_QQ_ASMPatchPCvsPCPATCH} and~\ref{fig:All_Timings_PP_ASMPatchPCvsPCPATCH}.  We see that the extended Vanka patches used here lead to higher times-to-solution at $k=8$, but that they remain comparable particularly with $\nu = 2$.  While we again have limited data to extrapolate from, we note that the time-to-solution now grows by a factor of about 12 from $k=4$ to $k=8$, which is still below the $\mathcal{O}(k^4)$ (or worse) that we might expect, but not quite as nice as seen for the Taylor-Hood cases.
%Fig12%%%%%%%%%%%%%%%%%%%%%%%%%%%%
%%%%%%%%%%%%%%%%%%%%%%%%%%%%%%%
%%%%%%%%%%%%%%%%%%%%%%%%%%%%%%%
%%%%%%%%%%%%%%%%%%%%%%%%%%%%%%%

%%%%%%%%%%%%%%%%%%%%%%%%%%%%%%%
%%%%%%%%%%%%%%%%%%%%%%%%%%%%%%%
\subsection{Stopping criteria}\label{ssec:stopping}

The results in the previous subsections indicate that the proposed monolithic multigrid methods are robust solvers for higher-order finite-element discretizations, particularly for the cases of Taylor-Hood discretizations on quadrilateral meshes and the $(\mathcal{BDM}_k,dP_{k-1})$ discretization on triangular grids.  In particular, we note that in both of these cases, the number of iterations needed for convergence to a fixed residual reduction tolerance by a relative factor of $10^{10}$ is relatively insensitive to both the polynomial order, $k$, and number of levels of refinement from a $5\times 5$ coarsest grid, denoted by $\ell$.  (See Figures~\ref{fig:All_QQ_ASMPatchPC_ItFGMRES} and~\ref{fig:All_BDMdP_ASMPatchPC_ItFGMRES}.)  We now tackle the question of how many iterations of preconditioned FGMRES are needed in order to realize the accuracy of these discretizations.  We focus on only these two discretizations here, but note that the same issue arises for all higher-order discretizations.

To do this, we consider the same model problem, with manufactured solution given in~\eqref{eq:model_problem}.  This problem has a smooth (trigonometric) velocity solution and zero pressure solution.  Consequently, we measure the quality of the velocity solution using the relative error in the ${\bf H^1}$ norm, while we measure the quality of the pressure solution using the absolute error in the $L^2$ norm.  \revise{For the velocity, we estimate the ${\bf H^1}$ norm of the difference between the computed solution and the exact solution using Firedrake's standard (adaptive) quadrature to evaluate the norms.  For the pressure, we orthogonalize the computed pressure against the constant function and similarly estimate the $L^2$ norm of the resulting function using Firedrake's standard (adaptive) quadrature.}

To find an appropriate stopping criterion for each value of $k$ and $\ell$, we solve the same linear system with increasingly stricter stopping tolerance on the relative residual norm, starting with the requirement that the residual norm be reduced below a factor of $10^{-2}$ times the original norm, and reducing this tolerance by successive factors of 2 until convergence in the measured errors is achieved.  We note this is somewhat imprecise, as the exact ``discretization error'' is unknown and varies depending on whether one focuses on only the quality of the velocity approximation (as is common practice) or also includes the pressure in consideration.  We also note this varies somewhat with the number of pre- and post-relaxation sweeps per cycle.  While we repeated this experiment with $1\leq \nu \leq 4$, we found that the optimal time-to-solution was typically achieved with $\nu = 1$, and there was only small variations in quality of solution with variations in $\nu$.

\begin{figure}
\graphicspath{{./Figures/}}
  \begin{center}
    \includegraphics[width=\textwidth]{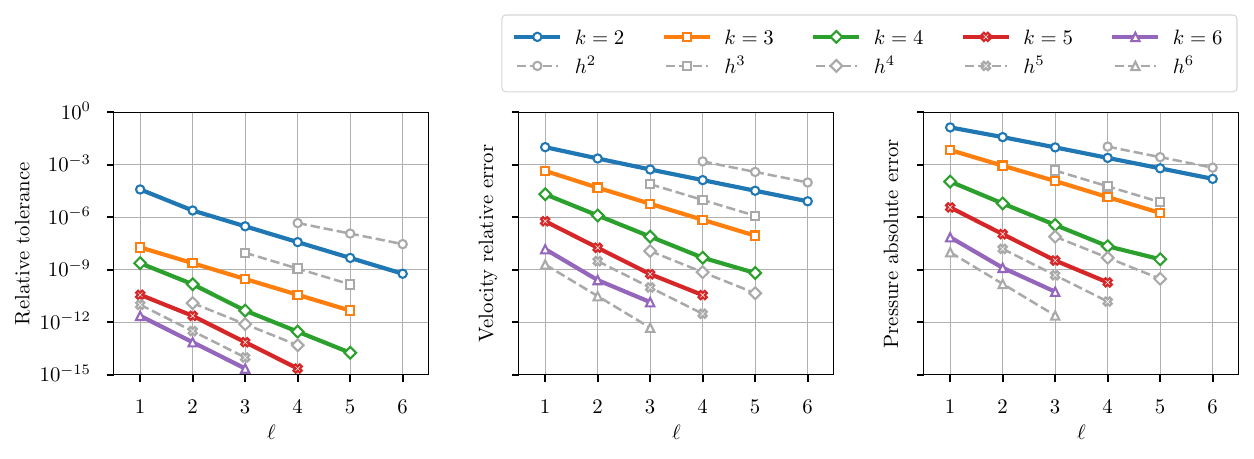}
    \caption{Selection of stopping criterion for the $(\mathcal{BDM}_k,dP_{k-1})$ discretization on triangular grids.  At {\bf left}, largest relative tolerance needed to achieve discretization error with varying order, $k$, and number of levels, $\ell$, in the multigrid hierarchy.  At {\bf center}, measured relative ${\bf H^1}$-norm errors in the velocity approximation at this stopping criterion.  At {\bf right}, measured absolute $L^2$-norm errors in the pressure approximation at this stopping criterion.}\label{fig:BDM_stopping}
  \end{center}
\end{figure}

Figure~\ref{fig:BDM_stopping} presents results of this study for the $(\mathcal{BDM}_k,dP_{k-1})$ discretization on triangular grids.  At left, we see a steady decrease in the relative tolerance required to achieve discretization-level accuracy as we increase $k$ and $\ell$.  Fitting the data, we see that for $k=2$ and $k=3$, a relative stopping criterion that scales like $h^3$ is sufficient to attain discretization error, but that this increases to $h^4$ for $k=4$ and $h^5$ for $k=5$ and $6$.  We note that we present data for all combinations of $k$ and $\ell$ that could be realized on our workstation, stopping experiments either when there was insufficient memory or the stopping tolerance needed to achieve discretization error fell below $10^{-15}$.  In the center and right panels of Figure~\ref{fig:BDM_stopping}, we show the measured errors in velocity and pressure approximations.  Standard finite-element convergence bounds for this discretization~\cite{JohnMixedCGDGStokes} predict convergence rates of $\mathcal{O}(h^k)$ for both velocity and pressure, which is consistent with our observations in this data.

\begin{figure}
\graphicspath{{./Figures/}}
  \begin{center}
    \includegraphics[width=\textwidth]{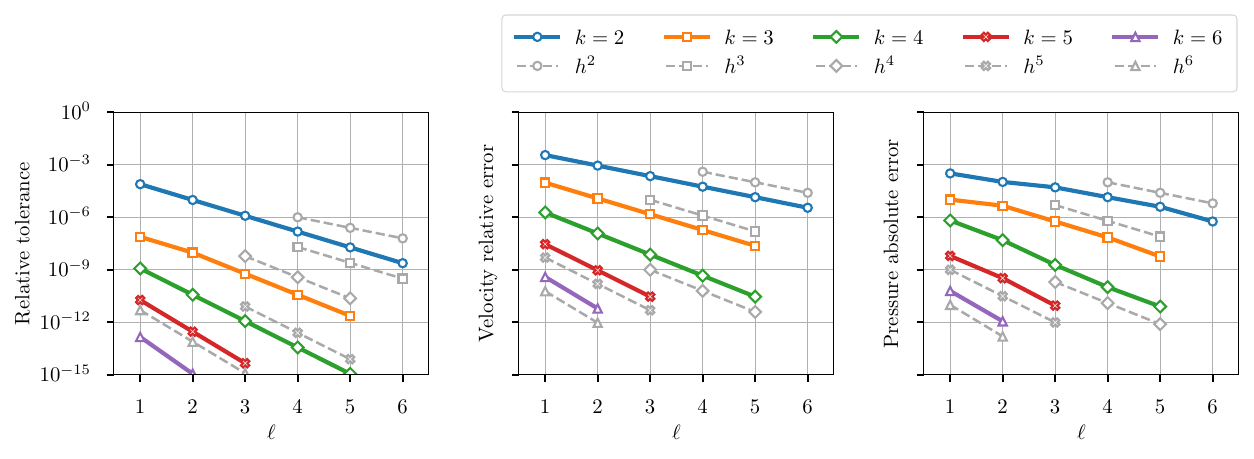}
    \caption{Selection of stopping criterion for the Taylor-Hood discretization on quadrilateral grids.  At {\bf left}, largest relative tolerance needed to achieve discretization error in the velocity approximation with varying order, $k$, and number of levels, $\ell$, in the multigrid hierarchy.  At {\bf center}, measured relative ${\bf H^1}$-norm errors in the velocity approximation at this stopping criterion.  At {\bf right}, measured absolute $L^2$-norm errors in the pressure approximation at this stopping criterion.}\label{fig:TH_stopping_velocity}
  \end{center}
\end{figure}

For the Taylor-Hood discretizations on quadrilaterals, we noticed a significant difference between the stopping criteria required to acheive convergence in the velocity approximation and those required to achieve convergence in the pressure approximation, so we present these results separately.  First, in Figure~\ref{fig:TH_stopping_velocity}, we consider the relative residual reduction tolerance needed to achieve discretization-level erorr in the velocity approximation.  Considering the left panel of the figure, we see that very similar tolerances are required for the lowest-order case $k=2$ as for the $(\mathcal{BDM}_k,dP_{k-1})$ discretization, but that stricter stopping tolerances appear necessary for larger $k$.  In fact, we see in general that our stopping tolerances scale like $\mathcal{O}(h^{k+1})$ for $2 \leq k \leq 6$.  With this, we observe $\mathcal{O}(h^k)$ convergence in both the velocity and pressure approximations, as expected from standard convergence theory~\cite{FEMIIGuermond}.

\begin{figure}
\graphicspath{{./Figures/}}
  \begin{center}
    \includegraphics[width=\textwidth]{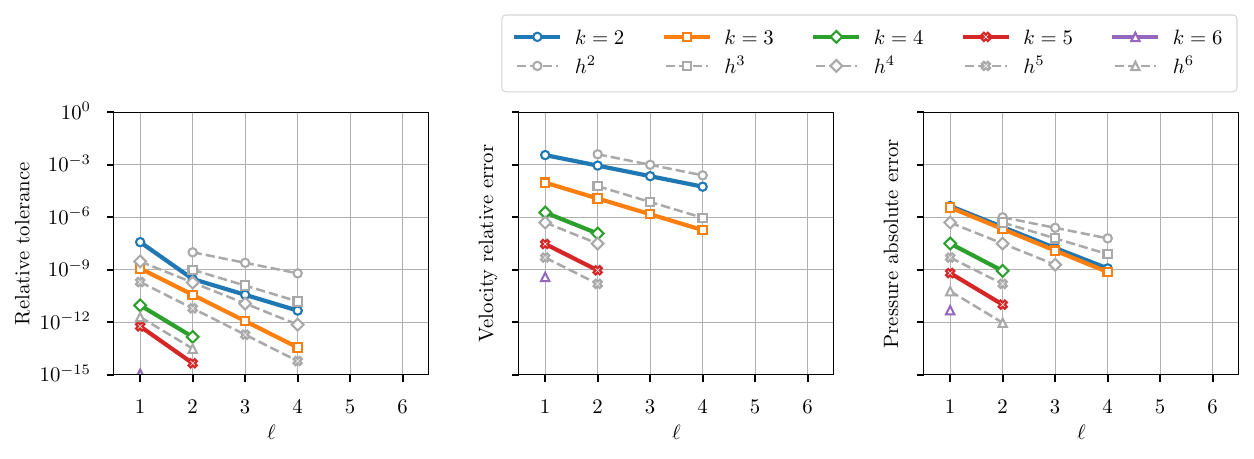}
    \caption{Selection of stopping criterion for the Taylor-Hood discretization on quadrilateral grids.  At {\bf left}, largest relative tolerance needed to achieve discretization error in the pressure approximation with varying order, $k$, and number of levels, $\ell$, in the multigrid hierarchy.  At {\bf center}, measured relative ${\bf H^1}$-norm errors in the velocity approximation at this stopping criterion.  At {\bf right}, measured absolute $L^2$-norm errors in the pressure approximation at this stopping criterion.}\label{fig:TH_stopping_pressure}
  \end{center}
\end{figure}

In preparing the numerical results presented in Figure~\ref{fig:TH_stopping_velocity}, we observed that the pressure error continued to be reduced as we further reduced the stopping criterion.  In many ways, this is not surprising, as we are approximating a zero pressure, which is exactly represented in the finite-element space.  In Figure~\ref{fig:TH_stopping_pressure}, we record the stopping criteria and corresponding measured errors when we choose the former to achieve full convergence of the pressure approximation.  We note, in particular, that these stopping criteria are much more strict than those above, scaling like $\mathcal{O}(h^3)$ for the lowest-order, $k=2$, case, but with much smaller constant, and scaling like $\mathcal{O}(h^{k+2})$ for $k \geq 3$.  The corresponding velocity errors do not change (as expected, since those in Figure~\ref{fig:TH_stopping_velocity} should represent the accuracy of the discretization).  However, we see notable improvement in the pressure errors in this experiment, scaling like $h^4$ for $k=2$ and $k=3$, and $h^{k+1}$ for $k > 3$.

Two important observations arise from these experiments.  First of all, the promise of higher-order discretizations does, indeed, pay off in this setting (with a smooth manufactured solution on a regular domain).  In Figure~\ref{fig:BDM_stopping}, for example, we see that the accuracy attained with $k=4$ on a $10\times 10$ mesh (with $\ell = 1$) is comparable to that attained with $k=2$ on a $160\times 160$ mesh (with $\ell = 5$).  This improvement in accuracy comes with a substantial reduction in time-to-solution, from about 45 seconds for the $k=2$ calculation on a large mesh, to about 2 seconds for $k=4$ on a coarse mesh.  Similar reductions in time-to-solution are seen consistently in these results, with the fastest times-to-solution for a given accuracy consistently seen for the highest-order discretizations that achieve that accuracy on the coarsest mesh.  Secondly, we note that the required stopping criteria depend dramatically on the polynomial order of the approximation, falling below what can be resolved in double-precision arithmetic for quite modest values of $k$ on grids with only a few tens of elements on a side.  This clearly signals that more work needs to be done in deriving robust stopping criteria for these discretizations, without which the practical utility of these discretizations is limited.

%%%%%%%%%%%%%%%%%%%%%%%
\section{Conclusion and future research direction}\label{sec:conclusions}

In this paper, we propose and study extensions to existing Vanka-style relaxation schemes for finite-element discretizations of the Stokes equations.  Following a topological extension of the ``pressure node oriented'' relaxation framework~\cite{JohnPatchesStokes}, we develop a composite Vanka relaxation scheme suitable for Taylor-Hood discretizations.  We also study the ``extended Vanka'' relaxation~\cite{ElementWiseExtendedVanka} for ${\bf H}(\text{div})$-$L^2$ discretizations, and demonstrate that both of these schemes lead to robust monolithic multigrid solvers with increasing polynomial order, with some dependence on the details of the discretization.  For the Taylor-Hood discretizations, we study two implementations of the relaxation and conclude that, even though matrix-full implementations have disadvantages for iterative solvers with weak preconditioners, having the fully assembled matrix leads to notable improvements in run time for these solvers over callback-based implementations.  Finally, we consider the issue of effective stopping criteria for these solvers, showing unfortunate dependence on the discretization order that clearly warrants further study.

Several other topics remain for future work.  A key question is whether the automated parameter choices considered here, based on estimating intervals for Chebyshev-accelerated relaxation, can be improved upon by local Fourier analysis or other techniques.  Also of interest is whether there are further improvements possible for solvers in this framework, either by improvements in implementation, improvements in efficiency of solving the patch systems, or alternative patch construction that leads to improvements in performance.  Finally, we note that many extensions of this work are possible, including to algebraic multigrid settings~\cite{AVoronin_etal_2023a} and to more complex systems of PDEs, such as Navier-Stokes or magnetohydrodynamics.

%%%%%%%%%%%%%%%%%%%%%%%%%%%%%%%%%%%%%%%%%%%%%%%%%%
%%%%%%%%%%%%%%%%%%%%%%%%%%%%%%%%%%%%%%%%%%%%%%%%%%
%\bibliographystyle{siamplain}
%\bibliographystyle(plain)

%\bibliographystyle{siamplain}
%\bibliography{references}

\section*{Acknowledgments}

This work was partially supported by an NSERC Discovery Grant. This work does not have any conflicts of interest.

\end{document}